\documentclass[a4paper,12pt]{article}
\usepackage{fancyhdr,graphicx}
\usepackage{amsfonts}
\usepackage{amssymb}
\usepackage{amsthm}
\usepackage{newlfont}
\usepackage{amsmath}
\usepackage[top=2.0cm,bottom=3.0cm,left=2.5cm,right=2.5cm]{geometry}
\usepackage{multirow}
\usepackage{array}
\usepackage{longtable}
\usepackage{bm}
\usepackage{latexsym}
\usepackage{amsmath}
\usepackage{graphicx}
\usepackage{amssymb}
\usepackage{mathrsfs}
\usepackage{setspace}
\usepackage{amssymb,amsfonts}
\allowdisplaybreaks \numberwithin{equation}{section}

\usepackage{lineno}
\usepackage{algorithm}
\usepackage{algorithmic}
\usepackage{xcolor}
\usepackage{float}
\usepackage{graphics}
\usepackage{graphicx}
\usepackage{epsfig}
\usepackage{booktabs,multirow}
\usepackage{diagbox}

\begin{document}

\title{{Associated Mersenne graphs\thanks{This work is supported by Shandong Provincial Natural Science Foundation (ZR2021MA071, ZR2019YQ02)
and National Natural Science Foundation of China (12171414).}}}

\author{
Jianxin Wei $^{a}$\footnote{Corresponding author.},
Yujun Yang $^{b}$\\
\scriptsize{$^{a}$
School of Mathematics and Statistics Science, Ludong University, Yantai, Shandong, 264025, P.R. China}\\
\scriptsize{$^{b}$
School of Mathematics and Information Science, Yantai University, Yantai, Shandong, 264005, P.R. China}\\
\scriptsize{E-mails:
wjx0426@ldu.edu.cn,
yangyj@ytu.edu.cn}}

\date{}

\maketitle

\begin{abstract}
In this paper,
we introduce the associated Mersenne graphs $\mathcal{M}_{n}$,
a new family of graphs determined by the definition of the Fibonacci-run graphs ({\"O}. E\v{g}ecio\v{g}lu, V. Ir\v{s}i\v{c}, 2021) in a circular manner.
The name of these graphs is justified with the fact that $|V(\mathcal{M}_{n})|$ is equal to the $n$-th associated Mersenne number.
Here,
we obtain several interesting structural and enumerative properties of associated Mersenne graphs including
the analogue of the fundamental recursion,
the number of vertices and edges,
the radius,
diameter,
center,
periphery and medianicity.
Finally,
several questions are
also listed.

\end{abstract}

\newcommand{\trou}{\vspace{1.5 mm}}
\newcommand{\noi}{\noindent}
\newcommand{\ol}{\overline}

\section{Introduction}
Let both $s$ and $t$ be integers and $s<t$.
Then we will use $s:t$ for the set $\{s,s+1,\ldots,t\}$ in this paper.
Let $B=\{0,1\}$ and $n\geq1$.
Then the elements from $\mathbf{B}_{n}=\{b_{1}b_{2}\ldots b_{n}\mid b_{i}\in B,i\in 1:n\}$ are called \emph{binary strings of length} $n$ (or simply strings).
The \emph{Hamming distance} $H(\alpha,$ $\beta)$ between two strings $\alpha$ and $\beta$ of $\mathbf{B}_{n}$ is the number of bits where $\alpha$ and $\beta$ differ.
The \emph{weight} of the string $\beta$ is $w(\beta)=\mathop{\sum}\limits^{n}_{t=1}b_{i}$.
We use the product notation for strings meaning concatenation,
for examples,
$1^{s}$ is the string of length $s$ with all bits equal $1$,
and $\mu\nu$ denotes the string by concatenating strings $\mu\in \mathbf{B}_{l}$ and $\nu\in\mathbf{B}_{m}$.
More generally,
$\mu^{s}$ is the concatenation of $s$ copies of $\mu$,
and if $\mathbf{S}$ is a subset of $\mathbf{B}_{n}$,
then $\mu \mathbf{S} \nu=\{\mu\gamma\nu|\gamma\in \mathbf{S}\}$.
Note that if $s=0$,
then $\mu^{s}=\lambda$,
where $\lambda$ is the null string.
A non-extendable sequence of contiguous equal bits in a string is called a \emph{run} (also a \emph{block}) of it,
for example,
$1^{3}0^{2}1^{2}0$ is a string of 4 runs.
A string $f$ is a \emph{factor} of a string $\alpha$
if $f$ appears as a sequence of $|f|$ consecutive bits of $\alpha$.

The hypercubes are one of the most popular architectures for interconnection networks for parallel computers.
From the view of graph theory,
the \emph{hypercube} $Q_{n}$ is the graph defined on the vertex set $\mathbf{B}_{n}$ and two vertices are adjacent when they have Hamming distance equal to $1$.
Clearly,
$Q_{n}$ has a regular degree $n$ and $V(Q_{n})=2^{n}$.
However,
the facts of high regularity and the exponentially growing number of vertices of $Q_{n}$,
as $n$ increases,
limit considerably the choice of the size of the networks.
These reasons led to an introduction of some alternative architectures such as Fibonacci
cube \cite{Hsu} and its various modifications.

A string is called a \emph{Fibonacci string} if it contains no two consecutive 1s.
Let $n\geq2$ and $\mathbf{F}_{n}$ be the set of all Fibonacci strings of length $n$,
i.e.,

$\mathbf{F}_{n}=\{b_{1}b_{2}\ldots b_{n}\in\mathbf{B}_{n}\mid b_{i}b_{i+1}=0,i\in1:n-1\}$.

\noindent
Then the \emph{Fibonacci cube} $\Gamma_{n}$ has $\mathbf{F}_{n}$ as the vertex set,
two vertices being adjacent if the Hamming distance of them equal to $1$.
For convenience,
it is always set $\mathbf{F}_{0}=\{\lambda\}$ and $\mathbf{F}_{1}=\{0,1\}$,
where $\lambda$ the null string.
So $\Gamma_{0}\cong K_{1}$ and $\Gamma_{1}\cong K_{2}$.
As mentioned above,
Fibonacci cubes were introduced as a model for interconnection networks.
However, they turned out to be interesting on their own,
see the survey of Klav\v{z}ar \cite{K1},
in which it can be found that Fibonacci cubes have several interesting properties on structure,
metric, combinatory and enumeration,
and are applicable in theoretical chemistry.

When Fibonacci cube was introduced,
it soon became increasingly popular,
and numerous variants and generalizations of Fibonacci cubes were proposed.
Let us just briefly mention a few of them,
such as
Fibonacci-run graph \cite{EI1,EI2},
$k$-Fibonacci cubes \cite{ESS2},
alternate Lucas cubes \cite{ESS1},
$k$-Lucas cubes \cite{ESS3},
generalized Fibonacci cube \cite{IKR1},
generalized Lucas cube \cite{IKR2},
daisy cube \cite{KM},
Lucas cube \cite{MPS} and
Pell graph \cite{Mu}.

Our principal motivation for the present paper comes from the Fibonacci-run graph and Lucas cube.
Lucas cubes is a class of graphs closely related to Fibonacci cubes.
A string is called a \emph{Lucas strings} if it is a Fibonacci string such that not both the first and the last bit are equal to $1$.
Let $n\geq1$ and $\mathbf{L}_{n}$ be the set of Lucas strings,
that is

$\mathbf{L}_{n}=\{b_{1}b_{2}\ldots b_{n}\in\mathbf{F}_{n}\mid b_{1}b_{n}=0\}$.

\noindent
Then the \emph{Lucas cube} $\Lambda_{n}$ is the graph defined on the vertex set $\mathbf{L}_{n}$,
and two vertices are adjacent when their Hamming distance equal to $1$.
Lucas cubes have been studied from different point of view \cite{CM,DTZ,IM,KM2}

The Fibonacci-run graph $\mathcal{R}_{n}$ has recently been introduced as an induced subgraph of hypercube by E\v{g}ecio\v{g}lu and Ir\v{s}i\v{c} \cite{EI1,EI2}.
A string of $\mathbf{B}_{n}$ is called \emph{run-constrained} if every run (also called a block) of 1s appearing in this string is immediately followed by a strictly longer run of 0s.
Let $\mathbf{R}_{n}$ be the set of all run-constrained strings of length $n$.
Then the graph $\mathcal{R}_{n}$ is defined on the vertex set $\mathbf{R}_{n}$,
and two vertices are adjacent when they have Hamming distance equal to $1$.
On graphs $\mathcal{R}_{n}$,
some basic properties were presented in \cite{EI1},
the nature of the degree sequences was studied in \cite{EI2},
and the radius and center were determined in \cite{Wei}.

Note that run-constrained strings of length $n\geq2$ must end with $00$.
Many times,
the vertices of Fibonacci-run graph $\mathcal{R}_{n}$ were viewed without the trailing pair of zeros in \cite{EI1,EI2}.
However,
for the convenience of the following study in the present paper,
the graph $\mathcal{R}_{n}$ is considered as a subgraph of $Q_{n}$ that the trailing 00 is always not omitted in the vertices of $\mathcal{R}_{n}$.
According to the definition of a run-constrained string,
it is natural to set $\mathbf{R}_{0}=\{\lambda\}$ and $\mathbf{R}_{1}=\{0\}$,
where $\lambda$ is the null string.

The graphs $\mathcal{R}_{n}$ (marked by the black dots) are shown in Fig. 1,
where $n\in1:6$.
Figures of a few large associated Mersenne graphs are given in the Appendix.

\begin{center}
\includegraphics[scale=0.90]{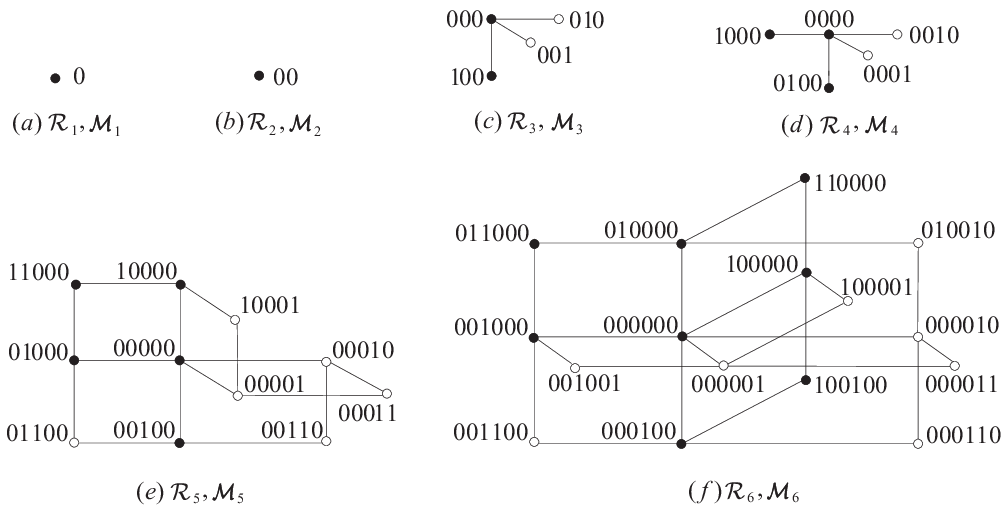}\\
{\footnotesize Fig. 1. The graphs $\mathcal{R}_{n}$ and $\mathcal{M}_{n}$, $n\in 1:6$.}
\end{center}

It is easy to see that a Fibonacci string of length $n$ is a string that contains no two consecutive
1s,
but a Lucas string is a string in which the substring $11$ is forbidden in a circular manner.
So the Fibonacci cube $\Gamma_{n}$ is the graph obtained from $Q_{n}$ by removing all vertices that contain $11$ as a factor,
but the Lucas cube $\Lambda_{n}$ is the graph obtained from $Q_{n}$ by removing all vertices which have a circulation containing $11$ as a factor.
Inspired by this,
this paper consider the circulation type of Fibonacci-run graph as follows.

A string of $\mathbf{B}_{n}$ is called \emph{run-constrained-circularly} if every run of 1s appearing in this string is immediately followed by a strictly longer run of 0s in a circular manner.
Let $\mathbf{M}_{n}$ be the set of all such words of length $n$.
Then the \emph{associated Mersenne graph} $\mathcal{M}_{n}$ is defined on the vertex set $\mathbf{M}_{n}$,
and two vertices adjacent if and only if their Hamming distance is 1.
For convenience,
we also set $\mathbf{M}_{0}=\{\lambda\}$ and $\mathbf{M}_{1}=\{0\}$.
Clearly,
a string of $\mathbf{R}_{n}$ must be a string of $\mathbf{M}_{n}$,
so the graph $\mathcal{R}_{n}$ is an induced subgraph of graph $\mathcal{M}_{n}$.

The name of these graphs $\mathcal{M}_{n}$ is justified with the fact that for any $n\geq1$,
$|V(\mathcal{M}_{n})|=|\mathbf{M}_{n}|=M_{n}$ (it will be shown in Section 3.1),
where $M_{n}$ are the \emph{associated Mersenne numbers} defined by the following recurrence for $n\geq2$
with the initial values $M_{0}=0,M_{1}=1$:

$M_{n}=M_{n-1}+M_{n-2}+1-(-1)^{n}$.

\noindent
These numbers form sequence A001350 in \cite{Sl} which will be introduced in Section 2.

The rest of the paper is organized as follows.
In Section 2,
some concepts and results needed in this paper are introduced.
In Section 3,
first,
the order of $\mathcal{M}_{n}$ is determined,
then,
a partition of $V(\mathcal{M}_{n})$ is given.
By this results,
the size of $\mathcal{M}_{n}$ is also determined and so the graph-theoretic decomposition of $\mathcal{M}_{n}$ is given.
In Section 4,
the radius, center, diameter and periphery of $\mathcal{M}_{n}$ are determined.
In Section 5,
the graphs $\mathcal{M}_{n}$ are proven been partial cubes and median graphs for only small values of $n$.
In Section 6,
several problems for further directions are listed.

\section{Preliminaries}

In this section some notion and results used in the paper are given.
Throughout the paper all graphs considered are finite and simple.
Other terminology and notation are introduced as they naturally occur in the following
and we use \cite{BM} for those not defined.

To avoid ambiguity with initial conditions,
Fibonacci numbers, Lucas numbers and associated Mersenne numbers are defined in the beginning of this section.

It is well known that Fibonacci numbers can be defined by the following recurrence formula with initial conditions
$F_{0}=0, F_{1}=1, F_{2}=1$:

$F_{n}=F_{n-1}+F_{n-2}$.

\noindent
The first values of Lucas numbers are 0, 1, 1, 2, 3, 5, 8, 13, 21, 34, 55, 89, 144, 233, 377, 610, 987, 1597, 2584.

For $n\geq2$,
Lucas numbers are defined as follows with initial conditions $L_{0}=2, L_{1}=1$ and $L_{2}=3$:

$L_{n}=L_{n-1}+L_{n-2}$.

\noindent
The first values of Lucas numbers are 2, 1, 3, 4, 7, 11, 18, 29, 47, 76, 123, 199, 322, 521, 843, 1364, 2207.

There is a relation often used between Fibonacci numbers and Lucas numbers.

\trou \noi {\bf Lemma 2.1 \cite{Sl}}.
$L_{n}=
\begin{cases}
(F_{n-k} + F_{n+k})/F_{k}, k>0~is~odd,\\
(F_{n+k} - F_{n-k})/F_{k}, k>0~is~even.\\
\end{cases}$

For example,
if $k=2$,
then $L_{n}=F_{n+2}-F_{n-2}$ and $L_{n-3}=F_{n-1}-F_{n-5}$.

As mentioned in the first section,
for $n\geq2$,
associated Mersenne numbers $M_{n}$ \cite{Sl} are defined as follows with initial conditions
$M_{0}=0, M_{1}=1$ and $M_{2}=1$:

$M_{n}=M_{n-1}+M_{n-2}+1-(-1)^{n}$.

\noindent
The first values of associated Mersenne numbers are 0, 1, 1, 4, 5, 11, 16, 29, 45, 76, 121, 199, 320, 521, 841, 1364, 2205.

It is easily seen that $M_{n}$ and $L_{n}$ satisfy the following relation.

\trou \noi {\bf Lemma 2.2 \cite{Sl}.}
\emph{For $n\geq0$.
$M_{n}=L_{n}-1-(-1)^{n}$.}

Next,
we turn to introduce some distance-based concept in graph theory.
Let $G$ be a (connected) graph.
Then the \emph{distance} $d_{G}(u,v)$ (or simply $d(u,v)$) between vertices $u$ and $v$ of $G$ is the length of a shortest path between
$u$ and $v$ in $G$.
It is well know that $d_{Q_{n}}(\alpha,\beta)=H(\alpha,\beta)$ for any two vertices $\alpha$ and $\beta$ of hypercube $Q_{n}$ \cite{HIK}.

The \emph{eccentricity} $e(v)$ of a vertex $v$ in a connected graph $G$ is the maximum distance between $v$ and another vertex, i.e.,

$e(v)=\max\limits_{u\in V(G)} d(v,u)$.

The \emph{diameter} $diam(G)$ of $G$ is the maximal eccentricity $e(v)$ when $v$ runs in $G$, i.e.,

$diam(G)=\max\limits_{v\in V(G)} e(v)=\max\limits_{u,v\in V(G)} d(u,v)$;

\noindent
A vertex $v$ is \emph{peripheral} when $e(v) = diam(G)$.
The \emph{periphery} $P(G)$ is the set of all peripheral vertices of $G$.

The \emph{radius} $rad(G)$ of $G$ is the minimum eccentricity of the vertices of $G$, i.e.,

$rad(G)=\min\limits_{v\in V(G)} e(v)$.

\noindent
A vertex is called \emph{central} if $e(v)=rad(G)$.
The \emph{center} $Z(G)$ of $G$ is the set of all central vertices.

The radius and center of graph $\mathcal{R}_{n}$ have been determined \cite{Wei},
but the diameter of $\mathcal{R}_{n}$ remains an open question \cite{EI2}.

In the last of this section,
some basic properties of Fibonacci-run graphs $\mathcal{R}_{n}$ are list.
Note that in paper \cite{EI1},
the trailing $00$ of every vertex of $\mathcal{R}_{n}$ is omitted,
the results introduced here have different appearance as they in \cite{EI1}.

For convenience,
it always set $|V(\mathcal{R}_{0})|=|\mathbf{R}_{0}|=1$.
For $n\geq1$,
the order of $\mathcal{R}_{n}$ is shown as follows.

\trou \noi {\bf Lemma 2.3\cite{EI1}.}
\emph{For $n\geq1$.
$|V(\mathcal{R}_{n})|(=|\mathbf{R}_{n}|)=F_{n}$.}

From Section 1,
it is known that the first few $V(\mathcal{R}_{n})(=\mathbf{R}_{n}$) are thus:


$\mathbf{R}_{1}=\{0\}$

$\mathbf{R}_{2}=\{00\}$

$\mathbf{R}_{3}=\{000,100\}$

$\mathbf{R}_{4}=\{0000,1000,0100\}$

$\mathbf{R}_{5}=\{00000,01000,00100,10000,11000\}$

$\mathbf{R}_{6}=\{000000,001000,000100,010000,100000,110000,011000,100100\}$.

\noindent
In general,
$\mathbf{R}_{n}$ can be got by the following lemma.

\trou \noi {\bf Lemma 2.4 \cite{EI1}.}
\emph{Let $n\geq1$.
Then the vertex set of $\mathcal{R}_{n}$ can be partitioned into}

$\mathbf{R}_{n}=\mathop{\bigcup}\limits^{\lceil n/2\rceil-1}_{i=0}1^{i}0^{i+1}\mathbf{R}_{n-2i-1}$

Further by considering the edges between subgraphs of $\mathcal{R}_{n}$ induced by $1^{i}0^{i+1}\mathbf{R}_{n-2i-1}$ ($i\in0:\lceil \frac{n}{2}\rceil-1$),
a schematic representation of the decomposition of the graph $\mathcal{R}_{n}$ is shown in Fig. 2,
where the arrowed line between two subsets of vertex denotes that there is an edge between every vertex of the beginning subset to the end subset.
\begin{center}
\includegraphics[scale=0.85]{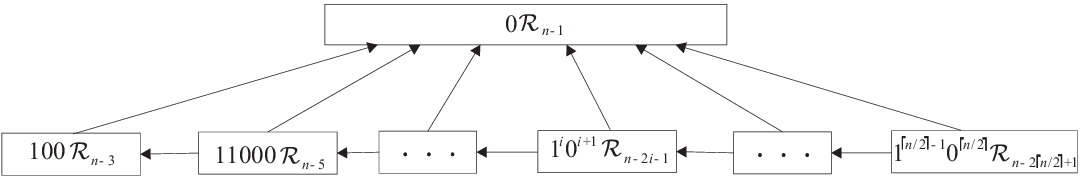}\\
{\footnotesize Fig. 2. The decomposition of the graph $\mathcal{R}_{n}$, $n=2k+1$ or $n=2k+2, k\geq1$.}
\end{center}

It is known that $|E(\mathcal{R}_{0})|=0,
|E(\mathcal{R}_{1})|=0,
|E(\mathcal{R}_{2})|=0,
|E(\mathcal{R}_{3})|=1,
|E(\mathcal{R}_{4})|=2.$
From Fig. 2,
the following calculation formula of $|E(\mathcal{R}_{n})|$ can be obtained.

\trou \noi {\bf Lemma 2.5\cite{EI1}.}
\emph{For $n\geq5$.
$|E(\mathcal{R}_{n})|=\mathop{\sum}\limits^{\lceil n/2\rceil-1}_{i=0}|E(\mathcal{R}_{n-2i-1})|
+|\mathbf{R}_{n-3}|
+2\mathop{\sum}\limits^{\lceil n/2\rceil-1}_{i=2}|\mathbf{R}_{n-2i-1}|$.}

Except the above recursion of $|E(\mathcal{R}_{n})|$,
there is a formula of $|E(\mathcal{R}_{n})|$ represented by Fibonacci numbers.

\trou \noi {\bf Lemma 2.6\cite{EI1}.}
\emph{If $n\geq6$,
then $|E(\mathcal{R}_{n})|=|E(\mathcal{R}_{n-1})|+|E(\mathcal{R}_{n-2})|+F_{n-3}+F_{n-5}$,
and if $n\geq8$,
then $|E(\mathcal{R}_{n})|=(3n-2)F_{n-8}+(5n-4)F_{n-7}$.}

\section{Decomposition of $\mathcal{M}_{n}$}

This section contains three subsections.
In the first subsection,
the order of graph $\mathcal{M}_{n}$ is determined.
In the second subsection,
a partition of vertex set is given for general $n$,
and so the vertex set of $\mathcal{M}_{n}$ can be obtained by the partition.
By counting the edges between the subsets of this partition,
the size of graph $\mathcal{M}_{n}$ is given in the third subsection.
Based on these results,
a decomposition of $\mathcal{M}_{n}$ is found.
It is useful for obtaining some interesting properties of $\mathcal{M}_{n}$.

\subsection{Order of $\mathcal{M}_{n}$}

The order of associated Mersenne graph $\mathcal{M}_{n}$ is determined for $n\geq0$ in this subsection.
Recall that a Lucas string is a string in which the factor $11$ is forbidden in a circular manner,
and $\mathbf{L}_{n}$ is the set of all Lucas strings of length $n$.

It is known that the first few $V(\mathcal{M}_{n})(=\mathbf{M}_{n})$ are thus:


$\mathbf{M}_{1}=\{0\}$,

$\mathbf{M}_{2}=\{00\}$,

$\mathbf{M}_{3}=\{000,001,010,100\}$,

$\mathbf{M}_{4}=\{0000,0001,0010,0100,1000\}$,

$\mathbf{M}_{5}=\{00000,00001,00010,00100,01000,10000,11000,01100,00110,00011,10001\}$.

\noindent
So $|V(\mathcal{M}_{0})|=0$,
$|V(\mathcal{M}_{1})|=1$,
$|V(\mathcal{M}_{2})|=1$,
$|V(\mathcal{M}_{3})|=4$,
$|V(\mathcal{M}_{4})|=5$, and
$|V(\mathcal{M}_{5})|=11$.
In general,
the order of $\mathcal{M}_{n}$ can be given as follows.

\trou \noi {\bf Theorem 3.1.}
\emph{Let $n\geq0$.
Then $|V(\mathcal{M}_{n})|=M_{n}$.}

\trou \noi {\bf Proof.}
Let $\nu$ be a vertex of graph $\mathcal{M}_{n}$.
Then $\nu$ can be circularly factored into the concatenation of factors $1^{i}0^{i+1}$ for some $i\geq0$.
Note that the factorization is unique,
and if $i=0$,
then $1^{i}0^{i+1}=0$.
For example,

$0001100010000111000011100=00)(0)(11000)(100)(0)(0)(1110000)(11100$

\noindent
Now for a factor $1^{i}0^{i+1}$,
we give a string $(10)^{i}0$ to correspond with it.
Thus,
for the vertex $\mu$,
a string $\mu'$ can be constructed.
For example,

$\mu=0001100010000111000011100$

$=00)(0)(11000)(100)(0)(0)(1110000)(11100$

$\rightarrow
\mu'=00)(0)(10100)(100)(0)(0)(1010100)(10101$

\noindent
It can be seen that the string $\mu'$ is a Lucas string.
Inspired by this,
We can give a bijection $\Phi$ between the vertex set of $\mathcal{M}_{n} (=\mathbf{M}_{n})$ and
a subset $\mathbf{L}'_{n}$ of the vertex set of Lucas cube $\Lambda_{n}(\mathbf{L}_{n})$,
where for $k\geq1$

$\mathbf{L}'_{n}=
\begin{cases}
\mathbf{L}_{n},n=2k+1\\
\mathbf{L}_{n}\setminus\{(10)^{k+1},(01)^{k+1}\},n=2k+2.\\
\end{cases}$

\noindent
It is clear that $\Phi$ is a bijection between $\mathbf{M}_{n}$ and $\mathbf{L}'_{n}$,
thus $|\mathbf{M}_{n}|=|\mathbf{L}'_{n}|$.
This means that

$|V(\mathcal{M}_{n})|=
\begin{cases}
L_{n}, n=2k+1,\\
L_{n}-2,n=2k+2\\
\end{cases}$

\noindent
thus $|V(\mathcal{M}_{n})|=L_{n}-1-(-1)^{n}=M_{n}$ by Lemma 2.2.
$\Box$

\subsection{Partition of $V(\mathcal{M}_{n})$}

A partition of $V(\mathcal{M}_{n})$ is given in this subsection.
This help us find the decomposition of graph $\mathcal{M}_{n}$.
For convenience,
for a string $\alpha=a_{1}a_{2}\ldots a_{t}$ of length $t$ and a string set $\mathbf{S}$,
let

$\overrightarrow{\alpha\mathbf{S}}=\{a_{i}\ldots a_{t}\mathbf{S}a_{1}\ldots a_{t-1}|i=1,2,\ldots,t\}$.

If $\mathbf{S}=\{\lambda\}$ (where $\lambda$ is the null string),
then $\overrightarrow{\alpha\mathbf{S}}$ is written as $\overrightarrow{\alpha}$ for brevity.
For examples,
$\overrightarrow{0\mathbf{S}}=0\mathbf{S}$,
$\overrightarrow{1110000\mathbf{S}}=1110000\mathbf{S}
\cup110000\mathbf{S}1\cup10000\mathbf{S}11\cup0000\mathbf{S}111
\cup000\mathbf{S}1110\cup00\mathbf{S}11100\cup0\mathbf{S}111000$.

\trou \noi {\bf Theorem 3.2.}
\emph{Let $n\geq1$.
Then $V(\mathcal{M}_{n})=\mathop{\bigcup}\limits^{\lceil n/2\rceil-1}_{i=0}\overrightarrow{1^{i}0^{i+1}\mathbf{R}_{n-2i-1}}$.}

\trou \noi {\bf Proof.}
Let $\nu=v_{1}v_{2}\ldots v_{n}$ be any vertex of $\mathcal{M}_{n}$.
Then the first bit of $\nu$ is either $0$ or $1$.
First we consider the case that $\nu$ starts with $0$.
Suppose,
in a circular manner,
that the length of the run containing $v_{1}(=0)$ is $i+1$ for some $i$,
where $0\leq i\leq\lceil\frac{n}{2}\rceil$,
and $v_{1}$ is the $(i-t+1)$-th bit of this run,
where $0\leq t\leq i$.
Then $\nu$ can be uniquely written as $00^{t}\omega1^{i}0^{i-t}$,
where $\omega$ is a run-constrained string of length $n-2i-1$.
This means that $\nu\in00^{t}\mathbf{R}_{n-2i-1}1^{i}0^{i-t}$.
Similarly,
for the case that $v_{1}=1$,
it can be shown that $\nu\in11^{t}0^{i+1}\mathbf{R}_{n-2i-1}1^{i-t-1}$,
where $1\leq i\leq \lceil\frac{n}{2}\rceil-1$ and $0\leq t<i$.
This yields the described vertex set decomposition.
$\Box$

\trou \noi {\bf Example 3.3.}
By Theorem 3.2,
we can obtain the vertex set of $\mathcal{M}_{n}$ for all $n\geq1$.
For $n=6$ and $7$,
$\mathbf{M}_{n}$ can be constructed as follows:

$\mathbf{M}_{6}=0\mathbf{R}_{5}\cup
(0\mathbf{R}_{3}10\cup00\mathbf{R}_{3}1\cup100\mathbf{R}_{3})
\cup(0\mathbf{R}_{1}1100\cup00\mathbf{R}_{1}110\cup000\mathbf{R}_{1}11
\cup1000\mathbf{R}_{1}1\cup11000\mathbf{R}_{1})$,

$\mathbf{M}_{7}=0\mathbf{R}_{6}\cup(0\mathbf{R}_{4}10\cup00\mathbf{R}_{4}1\cup100\mathbf{R}_{4})
\cup(0\mathbf{R}_{2}1100\cup00\mathbf{R}_{2}110\cup000\mathbf{R}_{2}11
\cup1000\mathbf{R}_{2}1\cup11000\mathbf{R}_{2})
\cup(0\mathbf{R}_{0}111000\cup00\mathbf{R}_{0}11100\cup000\mathbf{R}_{0}1110
\cup0000\mathbf{R}_{0}111\cup10000\mathbf{R}_{0}11\cup110000\mathbf{R}_{0}1
\cup1110000\mathbf{R}_{0})$,

\noindent
where $\mathbf{R}_{n}$ can be found in Section 2 for $n\in0:6$.

It is easily seen that $|\overrightarrow{1^{i}0^{i+1}\mathbf{R}_{n-2i-1}}|=(2i+1)|\mathbf{R}_{n-2i-1}|$.
So the following holds.

\trou \noi {\bf Corollary 3.4.}
\emph{Let $n\geq1$.
Then $|\mathbf{M}_{n}|=\mathop{\sum}\limits^{\lceil n/2\rceil-1}_{i=0}(2i+1)|\mathbf{R}_{n-2i-1}|$.}

This result implies that a associated Mersenne number can be represented as the sum of some Fibonacci numbers:
 $M_{n}=\mathop{\sum}\limits^{\lceil n/2\rceil-1}_{i=0}(2i+1)F_{n-2i-1}$ by Lemma 2.3.
To the best of our knowledge,
this equation was not studied earlier.

\subsection{Size of $\mathcal{M}_{n}$}

In this subsection,
we consider the size of $\mathcal{M}_{n}$.
From Theorem 3.2,
graph $\mathcal{M}_{n}$ contains the following $2i+1$ isomorphic subgraphs for every $i$ such that $0\leq i\leq \lceil\frac{n}{2}\rceil-1$:

$1^{i}0^{i+1}\mathbf{R}_{n-2i-1}$,
$\ldots$ ,
$10^{i+1}\mathbf{R}_{n-2i-1}1^{i-1}$,
$0^{i+1}\mathbf{R}_{n-2i-1}1^{i}$,
$\ldots$ ,
$0^{1}\mathbf{R}_{n-2i-1}1^{i}0^{i}$.

\noindent
This implies that there are $\mathop{\sum}\limits^{\lceil n/2\rceil-1}_{i=0}(2i+1)=(\lceil \frac{n}{2}\rceil)^{2}$ subgraphs of $\mathcal{M}_{n}$ to be considered.
For a given $i$,
inside each above subgraph the edges are inherited from the graph $\mathbf{R}_{n-2i-1}$.
Suppose that $\alpha$ and $\beta$ are two vertices contained in the above different subsets.
Then the $H(\alpha,\beta)\geq2$,
and so there is no edge between vertices of subgraphs induced by different vertex sets of $\overrightarrow{1^{i}0^{i+1}\mathbf{R}_{n-2i-1}}$.

Now we consider the edges between different subsets $\overrightarrow{1^{i}0^{i+1}\mathbf{R}_{n-2i-1}}$ and $\overrightarrow{1^{j}0^{j+1}\mathbf{R}_{n-2j-1}}$,
where $i\neq j$ and $i,j\in0:\lceil\frac{n}{2}\rceil-1$.
First the following result holds obviously.

\trou \noi {\bf Lemma 3.5.}
\emph{Let $1\leq i\leq \lceil\frac{n}{2}\rceil-1$.
Then all the strings of $\overrightarrow{1^{i-1}00^{i+1}\mathbf{R}_{n-2i-1}}$
and $\overrightarrow{01^{i-1}0^{i+1}\mathbf{R}_{n-2i-1}}$ are run-constrained-circularly strings of length $n$.}

Note that in Lemma 3.5, the sets $\overrightarrow{1^{i-1}00^{i+1}\mathbf{R}_{n-2i-1}}=\overrightarrow{01^{i-1}0^{i+1}\mathbf{R}_{n-2i-1}}$ holds if and only if $i=1$.
It can be seen that both the sets $\overrightarrow{1^{i-1}00^{i+1}\mathbf{R}_{n-2i-1}}$
and $\overrightarrow{01^{i-1}0^{i+1}\mathbf{R}_{n-2i-1}}$ can be obtained from $\overrightarrow{1^{i}0^{i+1}\mathbf{R}_{n-2i-1}}$ (by changing some one bit of the strings of $\overrightarrow{1^{i}0^{i+1}\mathbf{R}_{n-2i-1}}$ from 1 to 0).

It is known that two vertices $\alpha$ and $\beta$ of a associated Mersenne graph are adjacent if and only if $H(\alpha,\beta)=1$.
This means that if we change some one bit of vertex $\alpha$,
then obtain the vertex $\beta$.
It is easily seen that strings from $\overrightarrow{1^{i}0^{i+1}\mathbf{R}_{n-2i-1}}$ are of form $\overrightarrow{1^{i}0^{i+1}\omega}$,
where $\omega\in\mathbf{R}_{n-2i-1}$.
Every neighbor of a vertex of $\overrightarrow{1^{i}0^{i+1}\omega}$ can be obtained by changing some one bit of factor $1^{i}0^{i+1}$ or $\omega$.
Since all the neighbors such that obtained by changing some bit of $\omega$ are vertices of $\mathcal{R}_{n-2i-1}$,
we only need to consider the neighbors that obtained by changing some bit of the factor $1^{i}0^{i+1}$.
By Lemma 3.5,
if the first (resp. the $i$-th) bit of factor $1^{i}0^{i+1}$ is changed from $1$ to $0$,
then there is an edge between a vertex of $\overrightarrow{1^{i}0^{i+1}\omega}$ and a vertex of $\overrightarrow{1^{i-1}00^{i+1}\mathbf{R}_{n-2i-1}}$ (resp. $\overrightarrow{01^{i-1}0^{i+1}\mathbf{R}_{n-2i-1}}$).
It is obvious that if $i=1$,
then the first bit of factor $1^{i}0^{i+1}$ is also the $i$-th bit of it.
Note that both the string obtained by changing a bit of $0$ in $1^{i}0^{i+1}$ from $0$ to $1$,
and the string obtained by changing the $t$-th ($2\leq t\leq i-1$) bit of $1$ in $1^{i}0^{i+1}$ from $1$ to $0$,
are no longer a run-constrained-circularly string.
Thus,
the following result is got.

\trou \noi {\bf Lemma 3.6.}
\emph{Let $1\leq i\leq \lceil\frac{n}{2}\rceil-1$ and $0\leq j\leq \lceil\frac{n}{2}\rceil-1$.
If $i>1$,
then for a vertex $\nu$ of $\overrightarrow{1^{i}0^{i+1}\mathbf{R}_{n-2i-1}}$,
there are two neighbors of $\nu$ contained in some $\overrightarrow{1^{j}0^{j+1}\mathbf{R}_{n-2j-1}}$;
if $i=1$,
then for a vertex $\nu$ of $\overrightarrow{1^{i}0^{i+1}\mathbf{R}_{n-2i-1}}$,
there is only one neighbor of $\nu$ contained in some $\overrightarrow{1^{j}0^{j+1}\mathbf{R}_{n-2j-1}}$.}

By the above analysis,
the size of graph $\mathcal{M}_{n}$ can be get as follows.

\trou \noi {\bf Theorem 3.7.}
\emph{Let $n\geq5$.
Then}

$|E(\mathcal{M}_{n})|
=\mathop{\sum}\limits^{\lceil n/2\rceil-1}_{i=0}(2i+1)|E(\mathcal{R}_{n-2i-1})|
+3|\mathbf{R}_{n-3}|
+2\mathop{\sum}\limits^{\lceil n/2\rceil-1}_{i=2}(2i+1)|\mathbf{R}_{n-2i-1}|$.

By Lemma 2.5,
$|E(\mathcal{R}_{n-2i-1})|$ can be given.
Note that
for $n\in0:4$,
from Fig. 1 it is clear that
$|E(\mathcal{M}_{0})|=0,
|E(\mathcal{M}_{1})|=0,
|E(\mathcal{M}_{2})|=0,
|E(\mathcal{M}_{3})|=3,
|E(\mathcal{M}_{4})|=4.$
Combined with Theorem 3.7,
all the size of $\mathcal{M}_{n}$ can be determined.
There is a concise formula of the size of $\mathcal{M}_{n}$ as shown in Theorem 3.8.

To obtain the decomposition of graph $\mathcal{M}_{n}$,
it is needed to consider the explicit connecting relations between subgraphs induced by different $\overrightarrow{1^{i}0^{i+1}\mathbf{R}_{n-2i-1}}$.
Now for every $i (1\leq i\leq \lceil\frac{n}{2}\rceil-1)$,
we search the neighbors of a vertex of $\overrightarrow{1^{i}0^{i+1}\mathbf{R}_{n-2i-1}}$ according to the above analysis.
By Lemma 3.5,
it can be seen that for a given vertex $\nu$ of $\overrightarrow{1^{i}0^{i+1}\mathbf{R}_{n-2i-1}}$,
the set $\overrightarrow{1^{j}0^{j+1}\mathbf{R}_{n-2j-1}}$ which contain the neighbor of $\nu$ should satisfies with $j<i$.
By Lemma 3.6,
it should be suffice to find a neighbor of a vertex of $\overrightarrow{1^{i}0^{i+1}\mathbf{R}_{n-2i-1}}$ when $i=1$,
or two neighbors of a vertex of $\overrightarrow{1^{i}0^{i+1}\mathbf{R}_{n-2i-1}}$ for $i>1$  in some subset of $\overrightarrow{1^{j}0^{j+1}\mathbf{R}_{n-2j-1}}$ such that $j<i$.

First we consider the case $i=1$.
It can be seen that $\overrightarrow{000\mathbf{R}_{n-3}}\subseteq 0\mathbf{R}_{n-1}$.
For a vertex $\nu\in\overrightarrow{10^2\mathbf{R}_{n-3}}$,
there is a vertex $\nu'$ which can be obtained from $\nu$ by changing the bit $1$ in factor $10^2$ from $1$ to $0$.
It is obvious that $\nu'\in\overrightarrow{000\mathbf{R}_{n-3}}$ and $H(\nu,\nu')=1$.
This means that $\nu'\in0\mathbf{R}_{n-1}$,
and so there is a subgraph of $0\mathbf{R}_{n-1}$ isomorphic to the graph induced by $\overrightarrow{10^2\mathbf{R}_{n-3}}$.

Now we turn to consider the case $i>1$.
For convenience,
let $\overrightarrow{1^{i}0^{i+1}\mathbf{R}_{n-2i-1}}
=\mathbf{S_{1}}\cup\mathbf{S_{2}}\cup\mathbf{S_{3}}\cup\mathbf{S_{4}}$,
where
$\mathbf{S_{1}}=0\mathbf{R}_{n-2i-1}1^{i}0^{i}$,
$\mathbf{S_{2}}=00\mathbf{R}_{n-2i-1}1^{i}0^{i-1}$,
$\mathbf{S_{3}}=1^{i}0^{i+1}\mathbf{R}_{n-2i-1}$ and
$\mathbf{S_{4}}=\overrightarrow{1^{i}0^{i+1}\mathbf{R}_{n-2i-1}}\setminus(\mathbf{S_{1}}\cup \mathbf{S_{2}}\cup \mathbf{S_{3}})$.
It is easily seen $\mathbf{S_{4}}=0^{t}\mathbf{R}_{n-2i-1}1^{i}0^{i-t+1}\cup1^{s}0^{i+1}\mathbf{R}_{n-2i-1}1^{i-s}$,
where $3\leq t\leq i$ and $1\leq s\leq i-1$.

For $\mathbf{S_{1}}$,
since both $\mathbf{R}_{n-2i-1}1^{i-1}00^{i}$ and $\mathbf{R}_{n-2i-1}01^{i-1}0^{i}$ are run-constrained strings of length $n-1$ (that is, $\mathbf{R}_{n-2i-1}1^{i-1}00^{i}\subseteq \mathbf{R}_{n-1}$ and $\mathbf{R}_{n-2i-1}01^{i-1}0^{i}\subseteq \mathbf{R}_{n-1}$),
the set $(0\mathbf{R}_{n-2i-1}1^{i-1}00^{i}\cup0\mathbf{R}_{n-2i-1}01^{i-1}0^{i})\subseteq 0\mathbf{R}_{n-1}$.
So for every vertex $\nu\in \mathbf{S_{1}}$,
both the neighbors of $\mu$ contained in $0\mathbf{R}_{n-1}$.

For $\mathbf{S_{2}}$,
since $0\mathbf{R}_{n-2i-1}0\subseteq\mathbf{R}_{n-2i+1}$ and
$0\mathbf{R}_{n-2i-1}1^{i-1}00^{i-1}\subseteq \mathbf{R}_{n-1}$,
we know that the set $00\mathbf{R}_{n-2i-1}01^{i-1}0^{i-1}\subseteq 0\mathbf{R}_{n-2i+1}1^{i-1}0^{i-1}$ and $00\mathbf{R}_{n-2i-1}1^{i-1}00^{i-1}\subseteq 0\mathbf{R}_{n-1}$.
So for a vertex $\nu\in\mathbf{S_{2}}$,
one of the neighbors of $\nu$ belong to $0\mathbf{R}_{n-2i+1}1^{i-1}0^{i-1}(\subseteq\overrightarrow{1^{i-1}0^{i}\mathbf{R}_{n-2i+1}})$ and
the other to $0\mathbf{R}_{n-1}$.

For $\mathbf{S_{3}}$,
since the set $00\mathbf{R}_{n-2i-1}\subseteq\mathbf{R}_{n-2i+1}$ and $1^{i-1}0^{i+1}\mathbf{R}_{n-2i-1}\subseteq\mathbf{R}_{n-1}$,
it can be know that
$1^{i-1}00^{i+1}\mathbf{R}_{n-2i-1}\subseteq 1^{i-1}0^{i}\mathbf{R}_{n-2i+1}$ and $01^{i-1}0^{i+1}\mathbf{R}_{n-2i-1}\subseteq 0\mathbf{R}_{n-1}$.
So for a vertex $\nu\in\mathbf{S_{3}}$,
one of the neighbors of $\nu$ is in $1^{i-1}0^{i}\mathbf{R}_{n-2i+1}(\subseteq\overrightarrow{1^{i-1}0^{i}\mathbf{R}_{n-2i+1}})$ and the other in $0\mathbf{R}_{n-1}$.

For $\mathbf{S_{3}}$,
let $\mathbf{T_{1}}=0^{t}\mathbf{R}_{n-2i-1}1^{i}0^{i-t+1}(3\leq t\leq i)$ and $\mathbf{T_{2}}=1^{s}0^{i+1}\mathbf{R}_{n-2i-1}1^{i-s}(1\leq s\leq i-1)$.
Then for $\mathbf{T_{1}}$,
since $0^{t}\mathbf{R}_{n-2i-1}\subseteq 0^{t-2}\mathbf{R}_{n-2i+1}$ and $0^{t}\mathbf{R}_{n-2i-1}0\subseteq0^{t-1}\mathbf{R}_{n-2i+1}0$,
we know that the set
$0^{t}\mathbf{R}_{n-2i-1}1^{i-1}00^{i-t+1}\subseteq 0^{t-2}\mathbf{R}_{n-2i+1}1^{i-1}0^{i-t+2}$ and $0^{t}\mathbf{R}_{n-2i-1}01^{i-1}0^{i-t+1}\subseteq0^{t-1}\mathbf{R}_{n-2i+1}01^{i-1}0^{i-t+1}$.
So one of the neighbors of any vertex of $\mathbf{T_{1}}$ belongs to the set $0^{t-2}\mathbf{R}_{n-2i+1}1^{i-1}0^{i-t+2}(\subseteq\overrightarrow{1^{i-1}0^{i}\mathbf{R}_{n-2i+1}})$ and the other to $0^{t-1}\mathbf{R}_{n-2i+1}01^{i-1}0^{i-t+1}$$(\subseteq\overrightarrow{1^{i-1}0^{i}\mathbf{R}_{n-2i+1}})$.
For $\mathbf{T_{2}}$,
since $0^{2}\mathbf{R}_{n-2i-1}\subseteq\mathbf{R}_{n-2i+1}$
and $0\mathbf{R}_{n-2i-1}\subseteq\mathbf{R}_{n-2i+1}$,
it is know that the set $1^{s-1}00^{i+1}\mathbf{R}_{n-2i-1}1^{i-s}\subseteq1^{s-1}0^{i}\mathbf{R}_{n-2i+1}1^{i-s}$
and $1^{s}0^{i+1}\mathbf{R}_{n-2i-1}01^{i-s-1}\subseteq1^{s}0^{i}\mathbf{R}_{n-2i+1}01^{i-s-1}$.
This means that one of the neighbors of a vertex of $\mathbf{T_{2}}$ is in the set $1^{s-1}0^{i}\mathbf{R}_{n-2i+1}1^{i-s}\subseteq\overrightarrow{1^{i-1}0^{i}\mathbf{R}_{n-2i+1}}$ and the other in $1^{s}0^{i}\mathbf{R}_{n-2i+1}01^{i-s-1}\subseteq\overrightarrow{1^{i-1}0^{i}\mathbf{R}_{n-2i+1}}$.

By Theorem 3.2,
Lemmas 3.5 and 3.6 and the above analysis,
a schematic representation of the decomposition of graph $\mathcal{M}_{n}$ is shown in Fig. 3.
Compare the decompositions of $\mathcal{M}_{n}$ and $\mathcal{R}_{n}$ (as shown Fig. 2),
graph $\mathcal{R}_{n}$ is the subgraph of $\mathcal{M}_{n}$ induced by $\mathop{\bigcup}\limits^{\lceil n/2\rceil-1}_{i=0}1^{i}0^{i+1}\mathbf{R}_{n-2i-1}\subset V(\mathcal{M}_{n})$.
This also can be reflected in the contrast of Lemma 2.4 and Theorem 3.2,
Lemma 2.5 and Theorem 3.7.

\begin{center}
\includegraphics[scale=0.671]{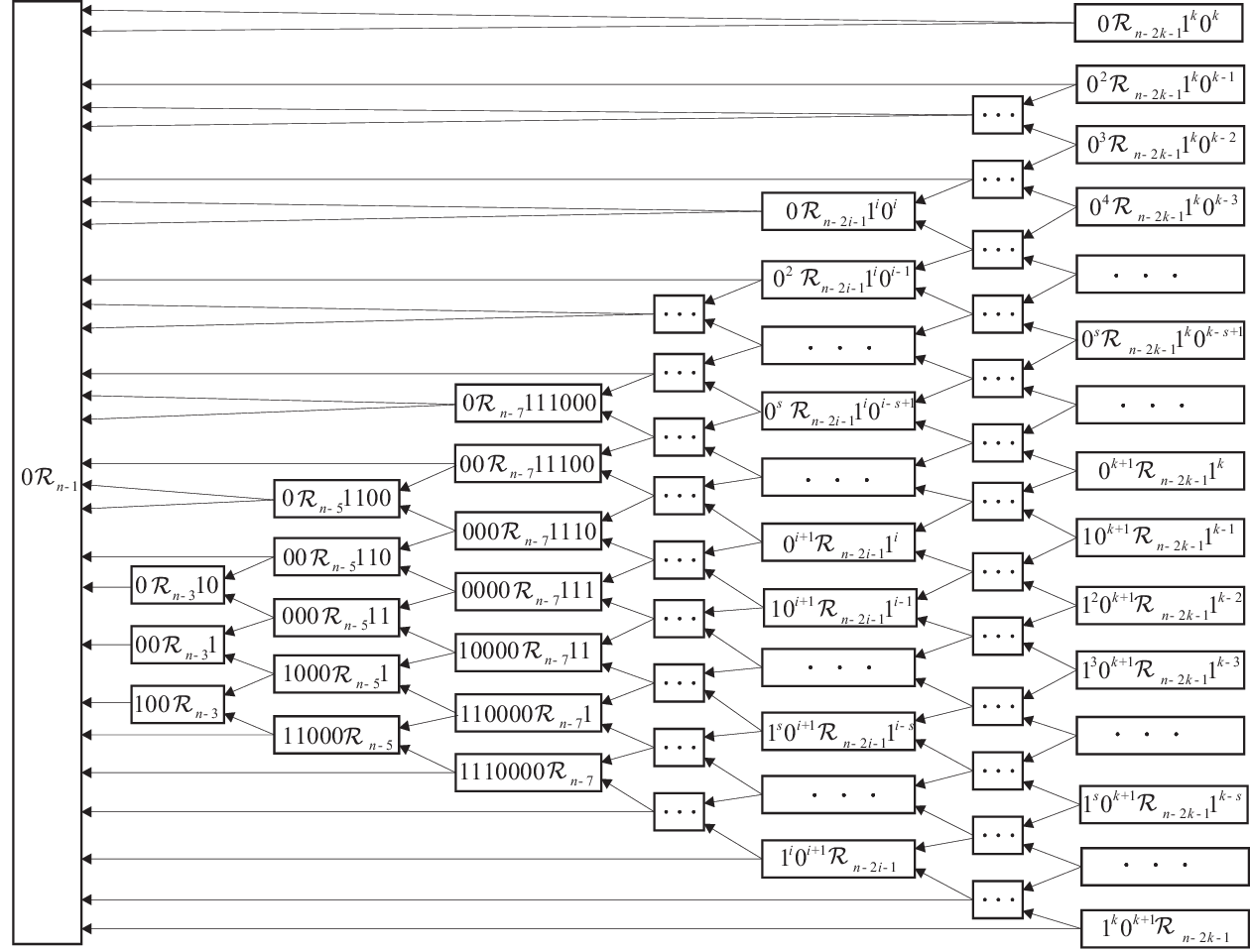}\\
{\footnotesize Fig. 3. The decomposition of the graph $\mathcal{M}_{n}$, $n=2k+1$ or $n=2k+2, k\geq1$.}
\end{center}

By using the relations among $F_{n}$, $L_{n}$ and $M_{n}$,
the formula of $|E(\mathcal{R}_{n})|$ and $|\mathbf{R}_{n}|$,
a much more simple formula of $|E(\mathcal{M}_{n})|$ can be shown as follows than its representation in Theorem 3.7.

\trou \noi {\bf Theorem 3.8.}
\emph{Let $n\geq4$.
Then $|E(\mathcal{M}_{n})|=nL_{n-3}$.}

\trou \noi {\bf Proof.}
The main method we use are induction.
It can be checked that $|E(\mathcal{M}_{4})|=4=4L_{1},
|E(\mathcal{M}_{5})|=15=5L_{2},
|E(\mathcal{M}_{6})|=24=6L_{3}$
and $|E(\mathcal{M}_{7})|=49=7L_{4}$.
So the equality holds for $n=4,5,6$ and $7$.
Now let $n\geq8$ and suppose that the equality $|E(\mathcal{M}_{n})|=nL_{n-3}$ holds for $n$.
Then by Theorem 3.7, Lemmas 2.5 and 2.6,
it follows that

$|E(\mathcal{M}_{n+2})|-|E(\mathcal{M}_{n})|$

$=|E(\mathcal{R}_{n+1})|+2|E(\mathcal{R}_{n})|+3F_{n-1}+5F_{n-3}$

$=|E(\mathcal{R}_{n+2})|+|E(\mathcal{R}_{n})|+2F_{n-1}+4F_{n-3}$

$=((3n+4)F_{n-6}+(5n+6)F_{n-5})+((3n-2)F_{n-8}+(5n-4)F_{n-7})+2F_{n-1}+4F_{n-3}$.

By the induction hypothesis and Lemma 2.1,
from the above equality,

$|E(\mathcal{M}_{n+2})|$

$=|E(\mathcal{M}_{n})|+((3n+4)F_{n-6}+(5n+6)F_{n-5})+((3n-2)F_{n-8}+(5n-4)F_{n-7})+2F_{n-1}+4F_{n-3}$

$=nL_{n-3}+((3n+4)F_{n-6}+(5n+6)F_{n-5})+((3n-2)F_{n-8}+(5n-4)F_{n-7})+2F_{n-1}+4F_{n-3}$

{\footnotesize $=n(F_{n-1}-F_{n-5})+((3n+4)F_{n-6}+(5n+6)F_{n-5})+((3n-2)F_{n-8}+(5n-4)F_{n-7})+2F_{n-1}+4F_{n-3}$}

=$(n+2)F_{n-1}+4F_{n-3}+(4n+6)F_{n-5}+(3n+4)F_{n-6}+(5n-4)F_{n-7}+(3n-2)F_{n-8}$

=$(n+2)F_{n-1}+4F_{n-3}+(4n+6)F_{n-5}+(6n+2)F_{n-6}+(2n-2)F_{n-7}$

=$(n+2)F_{n-1}+4F_{n-3}+(6n+4)F_{n-5}+(4n+4)F_{n-6}$

=$(n+2)F_{n-1}+4F_{n-3}+(4n+4)F_{n-4}+2nF_{n-5}$

=$(n+2)F_{n-1}+(2n+4)F_{n-3}+(2n+4)F_{n-4}$

=$(n+2)F_{n-1}+(2n+4)F_{n-2}$

=$(n+2)(F_{n-1}+F_{n-2})+(n+2)F_{n-2}$

=$(n+2)(F_{n}+F_{n-2})$

=$(n+2)L_{n-1}$,

\noindent
which completes the proof.
$\Box$

\trou \noi {\bf Corollary 3.9.}
\emph{Let $n\geq3$.
Then the generating function for $|E(\mathcal{M}_{n})|$ is}

$\mathop{\sum}\limits_{n\geq3}|E(\mathcal{M}_{n})|x^{n}
=\frac{(3-2x+4x^2-4x^3-3x^4)x^3}{(1-x-x^2)^2}$.

\trou \noi {\bf Proof.}
By Theorem 3.8,
$|E(\mathcal{M}_{n})|
=nL_{n-3}
=(n-3)L_{n-3}+3L_{n-3}$ for $n\geq4$.
Since $|E(\mathcal{M}_{3})|=3$ and the well-known fact that the generating function for $L_{n}$ is $\frac{2-x}{1-x-x^2}$,

$\mathop{\sum}\limits_{n\geq3}|E(\mathcal{M}_{n})|x^{n}$

$=3x^3+\mathop{\sum}\limits_{n\geq4}|E(\mathcal{M}_{n})|x^{n}$

$=3x^3+\mathop{\sum}\limits_{n\geq4}(n-3)L_{n-3}x^{n}+3\mathop{\sum}\limits_{n\geq4}L_{n-3}x^{n}$

$=3x^3+(\mathop{\sum}\limits_{n\geq4}(n-3)L_{n-3}x^{n-3}+3\mathop{\sum}\limits_{n\geq4}L_{n-3}x^{n-3})x^3$

$=3x^3+((\frac{2-x}{1-x-x^2}-2)'x+3(\frac{2-x}{1-x-x^2}-2))x^3$

$=\frac{(3-2x+4x^2-4x^3-3x^4)x^3}{(1-x-x^2)^2}$.
$\Box$

\section{Basic metric properties of $\mathcal{M}_{n}$}

By the decomposition of graph $\mathcal{M}_{n}$ as shown in Fig. 3,
we know that $\mathcal{M}_{n}$ is connected.
In this section,
$rad(\mathcal{M}_{n})$,
$Z(\mathcal{M}_{n})$,
$diam(\mathcal{M}_{n})$ and
$P(\mathcal{M}_{n})$ are determined.
We only consider the case that $\mathcal{M}_{n}$ is nontrivial graph,
that is,
$n\geq3$.

\trou \noi {\bf Lemma 4.1.}
\emph{Let $n\geq3$ and $\nu$ be a vertex of $\mathcal{M}_{n}$.
Then $e(\nu)\geq w(\nu)+\lceil\frac{n-w(\nu)}{2}\rceil-1$.}

\trou \noi {\bf Proof.}
Without loss of generality,
suppose $\nu\in V(\mathcal{M}_{n})$ is of form $\nu=1^{r_{1}}0^{s_{1}}\cdots1^{r_{t}}0^{s_{t}}$ by the circularity of a vertex of $\mathcal{M}_{n}$ for some $t\geq1$.

\noindent
Assume that among $s_{1},\ldots,s_{t}$,
there are $k$ odd numbers.
First we consider three special cases: $k=0,k=1$ and $k=2$.

\noindent
(1) $k=0$.
In this case,
$s_{i}$ is even for all $i\in1:t$.
Let

$\mu=0^{r_{1}}1^{\frac{s_{1}}{2}}0^{\frac{s_{1}}{2}}\cdots0^{r_{t}}1^{\frac{s_{t}}{2}}0^{\frac{s_{t}}{2}}$.

\noindent
Then $\mu$ is a vertex of $\mathcal{M}_{n}$ and $e(\nu)\geq d_{\mathcal{M}_{n}}(\nu,\mu)=w(\nu)+\lceil\frac{n-w(\nu)}{2}\rceil$.
For example,
if

$\nu=1001100001110000$, then

$\mu=0100011000001100$.

\noindent
(2) $k=1$.
Without loss of generality,
suppose $s_{t}$ is odd and $s_{i}$ is even for $i\in1:t\setminus\{t\}$.
Then let

$\mu=0^{r_{1}}1^{\frac{s_{1}}{2}}0^{\frac{s_{1}}{2}}
\cdots0^{r_{t}}1^{\lceil\frac{s_{t}}{2}\rceil-1}0^{\lceil\frac{s_{t}}{2}\rceil}$.

\noindent
It can be seen that $\mu$ is a vertex of $\mathcal{M}_{n}$ and $e(\nu)\geq d_{\mathcal{M}_{n}}(\nu,\mu)=w(\nu)+\lceil\frac{n-w(\nu)}{2}\rceil-1$.
For example,
if

$\nu=10011000011100000$, then

$\mu=01000110000011000$.

\noindent
(3) $k=2$.
Yet the general,
assume that $s_{1}$ and $s_{t}$ are odd and others are even.
Then let

$\mu=0^{r_{1}}1^{\lceil\frac{s_{t}}{2}\rceil}0^{\lceil\frac{s_{1}}{2}\rceil-1}
0^{r_{2}}01^{\frac{s_{2}}{2}}0^{\frac{s_{2}}{2}-1}
\cdots
0^{r_{t-1}}01^{\frac{s_{t-1}}{2}}0^{\frac{s_{t-1}}{2}-1}
0^{r_{t}}01^{\frac{s_{t}-1}{2}}0^{\frac{s_{t}-1}{2}}$.

\noindent
So $\mu$ is a vertex of $\mathcal{M}_{n}$ and $e(\nu)\geq d_{\mathcal{M}_{n}}(\nu,\mu)=w(\nu)+\lceil\frac{n-w(\nu)}{2}\rceil$.
For example,
if

$\nu=1100010011000011100000$, then

$\mu=0011000100011000001100$.

For convenience,
the strings of the form in cases (1), (2) and (3) are called strings of type (1), (2) and (3),
respectively.

Next we consider the most general cases, i.e.,
$k=2p+2$ is even and $k=2p+1$ is odd for $p\geq1$.
It is easily seen that if every vertex $\nu$ of $\mathcal{M}_{n}$ is shown as the concatenation of some of the above strings of types (1-3),
then the vertex $\mu$ can be constructed as the concatenation of the corresponding strings obtained in cases (1-3).
Thus,
the desired result holds.

In fact,
for $k=2p+2$,
suppose that $s_{i_{j}}$ is odd for $j\in1:2(p+1)$.
Then for $q\in1:p+1$,
the factor $1^{r_{i_{2q-1}}}0^{s_{i_{2q-1}}}\ldots1^{r_{i_{2q}}}0^{s_{i_{2q}}}$ of $\nu$ is of type (3),
and the factors segmented by the above factors are of type (1).
So $e(\nu)\geq w(\nu)+\lceil\frac{n-w(\nu)}{2}\rceil$.
For example,
note that the parts in italics or not are used to distinguish different factors here,
if

$\nu=100\emph{110001001000}1001110000\emph{1100010010011000}1110000$, then

$\mu=010\emph{001100010010}0100001100\emph{0011000101000010}0001100$.

\noindent
For $k=2p+1$,
suppose that $s_{i_{j}}$ is odd for $j\in1:2p+1$.
Then for $q\in1:p$,
the factor $1^{r_{i_{2q-1}}}0^{s_{i_{2q-1}}}\ldots1^{r_{i_{2q}}}0^{s_{i_{2q}}}$ of $\nu$ is of type (3),
the factor $1^{r_{i_{2p}+1}}0^{s_{i_{2p}}}\ldots1^{r_{i_{2p+1}}}0^{s_{i_{2p+1}}}$ of $\nu$ is of type (2)
and the factors segmented by the above factors are of type (1).
So $e(\nu)\geq w(\nu)+\lceil\frac{n-w(\nu)}{2}\rceil-1$.
For example,
if

$\nu=100\emph{110001001000}1001110000\emph{1100010010011000}111000011000\emph{100}$, then

$\mu=010\emph{001100010010}0100001100\emph{0011000101000010}000110000100\emph{010}$.

\noindent
This completes the proof.
$\Box$

\trou \noi {\bf Lemma 4.2.}
\emph{Let $n\geq3$,
$\alpha$ and $\beta$ be vertices of $\mathcal{M}_{n}$.
Then $d_{\mathcal{M}_{n}}(\alpha,\beta)\leq w(\alpha)+w(\beta)$.}

\trou \noi {\bf Proof.}
Suppose $\alpha=a_{1}a_{2}\ldots a_{n} \in V(\mathcal{M}_{n})$ and $a_{i_{1}}=a_{i_{2}}=\cdots=a_{i_{w}}=1$,
where $w(>0)$ is the weight of $\alpha$ and $i_{1}\leq i_{2}\leq\cdots\leq i_{w}$.
Let $j\in 1:w$ and $\alpha_{j}$ be the string obtained by changing the $i_{1},\ldots,i_{j}$-th bits of $\alpha$ from $1$ to $0$.
Then for $j\in 1:w$,
$\alpha_{j}\in V(\mathcal{M}_{n})$ and $\alpha_{w}=0^{n}$.
So there is a $\alpha,0^{n}$-path of length $w$:
$\alpha-\alpha_{1}-\alpha_{2}-\ldots-\alpha_{w}$.
Similarly,
for any vertex $\beta\in V(\mathcal{M}_{n})$,
a $\beta,0^{n}$-path of length $w(\beta)$ can be obtained.
This means that for any vertices $\alpha,\beta\in V(\mathcal{M}_{n})$,
there is a $\alpha,\beta$-path of length $w(\alpha)+w(\beta)$.
$\Box$

From the above proof we know that for every vertex of $\mathcal{M}_{n}$ with weight more that $0$,
there is a path connecting the vertex with $0^{n}$.
So we know that $\mathcal{M}_{n}$ is connected.
As mentioned above,
this is also indicated from the decomposition of $\mathcal{M}_{n}$ shown in Fig. 3.

\trou \noi {\bf Lemma 4.3.}
\emph{Let $n\geq3$ and $\nu\in V(\mathcal{M}_{n})$.
If $w(\nu)=0$, then $e(\nu)=\lceil\frac{n}{2}\rceil-1$,
if $w(\nu)=1$, then $e(\nu)=\lceil\frac{n}{2}\rceil$.}

\trou \noi {\bf Proof.}
If $w(\nu)=0$,
then $\nu=0^{n}$.
It is obvious that there are at most $\lceil\frac{n}{2}\rceil-1$ ones contained in a vertex of $\mathcal{M}_{n}$,
one of such vertex can be chosen as $\mu=1^{\lceil\frac{n}{2}\rceil-1}0^{\lceil\frac{n}{2}\rceil}$.
By Lemma 4.2,
it is known that $e(\nu)\leq \lceil\frac{n}{2}\rceil-1$.
It is easily seen that $e(\nu)\geq d_{\mathcal{M}_{n}}(\nu,\mu)=\lceil\frac{n}{2}\rceil-1$.
So $e(\nu)=\lceil\frac{n}{2}\rceil-1$.

Now we suppose $w(\nu)=1$.
Then $\nu$ can be taken as $0^{t}10^{n-t-1}$ for $t\in0:n-1$.
For every $t$,
a vertex $\mu$ of $\mathcal{M}_{n}$ can be constructed as follows.
If $t\in0:\lceil\frac{n}{2}\rceil-1$,
then let $\mu=1^{t}0^{n-\lceil\frac{n}{2}\rceil+1}1^{\lceil\frac{n}{2}\rceil-t-1}$.
If $t\in \lceil\frac{n}{2}\rceil:n-1$,
then let $\mu=0^{t-\lceil\frac{n}{2}\rceil+1}1^{\lceil\frac{n}{2}\rceil-1}0^{n-t}$.
It is obvious that $e(\nu)\geq d_{\mathcal{M}_{n}}(\nu,\mu)=\lceil\frac{n}{2}\rceil$.
Since the string $\mu$ is of the most number of $1$s and $w(\mu)=\lceil\frac{n}{2}\rceil-1$,
$e(\nu)\leq w(\nu)+w(\mu)=1+\lceil\frac{n}{2}\rceil-1=\lceil\frac{n}{2}\rceil$ by Lemma 4.2.
Thus,
we know $e(\nu)=d_{\mathcal{M}_{n}}(\nu,\mu)=\lceil\frac{n}{2}\rceil$.
$\Box$

\trou \noi {\bf Theorem 4.4.}
\emph{Let $n\geq3$.
Then $rad(\mathcal{M}_{n})=\lceil\frac{n}{2}\rceil-1$,
and $Z(\mathcal{M}_{n})=\{0^{n}\}$.}

\trou \noi {\bf Proof.}
By Lemma 4.1,
it is known for any vertex $\nu\in V(\mathcal{M}_{n})$ of weight $w$,
$e(\nu)\geq w+\lceil\frac{n-w}{2}\rceil-1$.
By Lemma 4.3,
$e(0^{n})=\lceil\frac{n}{2}\rceil-1$.
So $0^{n}\in Z(\mathcal{M}_{n})$.
It is easily seen that for a given $n$,
$f(w)=w+\lceil\frac{n-w}{2}\rceil-1$ is a monotone increasing function of $w(\geq0)$.
Furthermore,
if $n$ is even,
then $f(w)>f(0)$ for $w\geq1$;
if $n$ is odd,
then $f(w)>f(1)=f(0)$ for $w\geq2$.
By Lemma 4.3,
it is known that if $w(\nu)=1$,
then $e(\nu)=\lceil\frac{n}{2}\rceil>f(1)$.
Thus,
there is only one vertex $0^{n}\in Z(\mathcal{M}_{n})$.
This completes the proof.
$\Box$

\trou \noi {\bf Theorem 4.5.}
\emph{Let $n\geq3$.
Then $diam(\mathcal{M}_{n})=2(\lceil\frac{n}{2}\rceil-1)$.
If $n$ is odd,
then $P(\mathcal{M}_{n})=\overrightarrow{1^{\lceil\frac{n}{2}\rceil-1}0^{\lceil\frac{n}{2}\rceil}}$.
If $n$ is even,
then for $t\geq1$,}

$P(\mathcal{M}_{n})=
\begin{cases}
\overrightarrow{10^{2}1^{2t-2}0^{2t-1}}
\cup\overrightarrow{1^{t-1}0^{t}1^{t}0^{t+1}}
\cup\overrightarrow{1^{2t-1}0^{2t+1}},~n=4t\\
\overrightarrow{10^{2}1^{2t-1}0^{2t}}
\cup\overrightarrow{1^{t-1}0^{t}1^{t+1}0^{t+2}}
\cup\overrightarrow{1^{t}0^{t+1}1^{t}0^{t+1}}
\cup\overrightarrow{1^{2t}0^{2t+2}},~n=4t+2.\\
\end{cases}$

\trou \noi {\bf Proof.}
It is clear that the maximum weight of the vertex of $\mathcal{M}_{n}$ is $\lceil\frac{n}{2}\rceil-1$.
Since $d_{\mathcal{M}_{n}}(\alpha,\beta)\leq w(\alpha)+w(\beta)\leq2(\lceil\frac{n}{2}\rceil-1)$ for any vertices $\alpha,\beta\in V(\mathcal{M}_{n})$ by Lemma 4.2,
it is known that $diam(\mathcal{M}_{n})\leq2(\lceil\frac{n}{2}\rceil-1)$.

It is only need to show that there are vertices $\alpha$ and $\beta$ such that $d_{\mathcal{M}_{n}}(\alpha,\beta)=2(\lceil\frac{n}{2}\rceil-1)$.
In fact,
let $\alpha=1^{\lceil\frac{n}{2}\rceil-1}0^{\lceil\frac{n}{2}\rceil}$ and $\beta=0^{\lceil\frac{n}{2}\rceil}1^{\lceil\frac{n}{2}\rceil-1}$.
Then $d_{\mathcal{M}_{n}}(\alpha,\beta)=2(\lceil\frac{n}{2}\rceil-1)$.
This means that $diam(\mathcal{M}_{n})=2(\lceil\frac{n}{2}\rceil-1)$.

A vertex $\alpha$ of $\mathcal{M}_{n}$ is peripheral only if it is satisfied with that
$w(\alpha)=\lceil\frac{n}{2}\rceil-1$ and
there is a vertex $\beta$ such that $d_{\mathcal{M}_{n}}(\alpha,\beta)=2(\lceil\frac{n}{2}\rceil-1)$.
It is distinguished two cases to find $P(\mathcal{M}_{n})$ by $n$ is odd or even.

First we consider the case $n=2k+1$, $k\geq1$.
Suppose that $\alpha$ is a peripheral vertex.
Then $w(\alpha)=\lceil\frac{n}{2}\rceil-1=k$,
and so $\alpha$ has exactly two runs in a circular manner.
Hence we know that $\alpha$ is of form $0^{t}1^{k}0^{k+1-t}$ for $t\in1:k+1$ or $1^{s}0^{k+1}1^{k-s}$ for $s\in1:k$,
that is,
$\alpha\in \overrightarrow{1^{k}0^{k+1}}$.
For $\alpha=0^{t}1^{k}0^{k+1-t}$,
let $\beta=1^{t-1}0^{k+1}1^{k-t+1}$;
for $\alpha=1^{s}0^{k+1}1^{k-s}$,
let $\beta=0^{s}1^{k}0^{k-s+1}$.
Then $d_{\mathcal{M}_{n}}(\alpha,\beta)=2k=2(\lceil\frac{n}{2}\rceil-1)$.
Thus,
$P(\mathcal{M}_{n})=\overrightarrow{1^{k}0^{k+1}}$ if $n$ is odd.

Now we turn to consider the case $n=2k+2$ for $k\geq1$.
In this case,
$diam(\mathcal{M}_{n})=2k=n-2$.
Suppose that $\alpha$ is a peripheral vertex of $\mathcal{M}_{n}$.
Then $w(\alpha)=k$ and so $\alpha$ has either two runs or four runs in a circular manner.

By using similar method to the case $n$ being odd,
it can be shown that if $\alpha$ has exactly two runs,
then $\alpha\in \overrightarrow{1^{k}0^{k+2}} \subseteq P(\mathcal{M}_{n})$.

Suppose $\alpha$ has four runs and so it is of the form $1^{s}0^{s+1}1^{p}0^{p+1}$ for some $s,p\geq1$.
Let $\beta\in V(\mathcal{M}_{n})$ and $d_{\mathcal{M}_{n}}(\alpha,\beta)=2k=n-2$.
Then $\beta\in \overrightarrow{1^{m}0^{m+1}1^{l}0^{l+1}}$ for some $m\geq0$ and $l\geq1$,
and there are exactly two coordinates $i$ and $j$ in which $\alpha$ and $\beta$ have the same bits.
Without loss of generality,
assume $j=i+1$ or $j>i+1$.
According to the positions of $i$ and $j$,
we can find the vertices $\alpha$ and $\beta$ by distinguishing $m>1$ and $m=0$ as follows.

If $m>1$,
then for $j=i+1$,
there are two possibilities of $i$ by the circularity of the vertex of $\mathcal{M}_{n}$:
$i=s+1$ or $i=2s+2p+2$.
For $i=s+1$,
we have

$\alpha=1^{s}\emph{00}0^{s-1}1^{p}0^{p+1}$, and

$\beta=0^{l-1}\emph{00}1^{m}0^{m+1}1^{l}$,

\noindent
where the parts in italics denotes the coordinates $i$ and $j$.
So we have $s=l-1,s+1=m,p=m+1$ and $p+1=l$.
From those equalities,
$s=p$ and $l=m+2$.
Let $s=t$.
Then $p=t,l=t+1$ and $m=t-1$.
So
$\alpha=1^{t}0^{t+1}1^{t}0^{t+1}$ and
$\beta=0^{t+2}1^{t-1}0^{t}1^{t+1}$.
For $i=2s+2p+2$,
we have

$\alpha=1^{s}0^{s-1}\emph{00}1^{p}0^{p+1}$, and

$\beta=0^{l+1}1^{m}\emph{00}0^{m-1}1^{l}$,

\noindent
So $s=l+1,s-1=m,p=m-1$ and $p+1=l$.
From those equalities,
$s=p+2$ and $l=m$.
Let $l=t$.
Then $m=t,s=t+1$ and $p=t-1$.
Hence,
$\alpha=1^{t+1}0^{t+2}1^{t-1}0^{t}$ and
$\beta=0^{t+1}1^{t}0^{t+1}1^{t}$.

If $m>1$ and $j>i+1$,
then there are three possibilities of $i$ and $j$ by the circularity of the vertex of $\mathcal{M}_{n}$:
$i=s+1$ and $j=2s+2$;
$i=s+1$ and $j=2s+p+2$; and
$i=s+1$ and $j=2s+2p+2$.
For the first case,
we have

$\alpha=1^{s}\emph{0}0^{s-1}\emph{0}1^{p}0^{p+1}$, and

$\beta=0^{l}\emph{0}1^{m}\emph{0}0^{m}1^{l}$.

\noindent
So $s=l,s-1=m,p=m$ and $p+1=l$.
By those equalities,
we get $s=p+1$ and $l=m+1$.
Let $s=t$.
Then $p=t-1,m=t-1$ and $l=t$.
Hence,
$\alpha=1^{t}0^{t+1}1^{t-1}0^{t}$, and
$\beta=0^{t+1}1^{t-1}0^{t}1^{t}$.
For the second case, by

$\alpha=1^{s}\emph{0}0^{s}1^{p}\emph{0}0^{p}$, and

$\beta=0^{l}\emph{0}1^{m}0^{m}\emph{0}1^{l}$,

\noindent
we have $s=p$ and $l=m$.
Let $s=t$.
Then $\alpha=1^{t}0^{t+1}1^{t}0^{t+1}$ and
$\beta=0^{t+1}1^{t}0^{t+1}1^{t}$.
For the third case,

$\alpha=1^{s}\emph{0}0^{s}1^{p}0^{p}\emph{0}$ and

$\beta=0^{l-1}\emph{0}1^{m}0^{m+1}1^{l}\emph{0}$.

\noindent
So $p=s+1$ and $l=m+1$ can be got.
Let $p=t$.
Then $\alpha=1^{t-1}0^{t}1^{t}0^{t+1}$ and
$\beta=0^{t}1^{t-1}0^{t}1^{t}0$.

Similarly,
if $m=0$,
then only for the case $s=1$, $i=2s$ and $j=i+1$,
there are vertices
$\alpha=10^{2}1^{k-1}0^{k}$ and
$\beta=00^{2}0^{k-1}1^{k}$.

Note that if $\alpha$ is a peripheral vertex,
then $\overrightarrow{\alpha}\subseteq P(\mathcal{M}_{n})$.
From the above discussions,
we know that for $n=2k+2$ ($k\geq1$),
$\overrightarrow{1^{k}0^{k+2}}\subseteq P(\mathcal{M}_{n})$,
$\overrightarrow{10^{2}1^{k-1}0^{k}}\subseteq P(\mathcal{M}_{n})$.
Further,
let $k=2t-1(t\geq1)$.
Then for $n=2k+2=4t$,
$\overrightarrow{1^{t-1}0^{t}1^{t}0^{t+1}}\subseteq P(\mathcal{M}_{n})$.
Let $k=2t(t\geq1)$.
Then for $n=2k+2=4t+2$,
$\overrightarrow{1^{t-1}0^{t}1^{t+1}0^{t+2}}
\cup\overrightarrow{1^{t}0^{t+1}1^{t}0^{t+1}}\subseteq P(\mathcal{M}_{n})$.
This complete the proof.
$\Box$

\section{$\mathcal{M}_{n}$ as partial cubes and median graphs}

Let $H$ and $G$ be connected graphs.
Then a mapping $f: V(H)\rightarrow V(G)$ is an \emph{isometric embedding} if $d_{H}(u,v)=d_{G}(f(u),f(v))$ for any $u,v\in V(H)$.
A \emph{partial cube} is a connected graph that admits an isometric embedding into a hypercube \cite{HIK}.

A \emph{median} of vertices $u, v,w \in V(G)$ is a vertex of $G$ that simultaneously lies
on a shortest $u, v$-path, a shortest $u,w$-path, and a shortest $v,w$-path.
The graph $G$ is called a \emph{median graph} if every triple of its vertices has a unique median.
It is known that a median graph must be a partial cube,
and hypercube $Q_{n}$ is a median graph for every $n\geq1$ \cite{HIK}.

In this section,
it is shown graphs $\mathcal{M}_{n}$ are in general not partial cubes,
and so they are not median graphs.

\trou \noi {\bf Theorem 5.1.}
\emph{Graph $\mathcal{M}_{n}$ is a partial cube if and only if $n\leq8$.}

\trou \noi {\bf Proof.}
Suppose $\mathcal{M}_{n}$ is a partial cube.
Then for $\alpha,\beta\in V(\mathcal{M}_{n})$,
$d_{\mathcal{M}_{n}}(\alpha,\beta)=H(\alpha,\beta)$.
For $n\geq9$,
there exist vertices $\mu=1111000000^{n-9}$ and $\nu=1001000000^{n-9}$.
It is easily seen that $H(\mu,\nu)=2$.
By contradiction suppose that $\mathcal{M}_{n}$ is a partial cube for $n\geq9$,
then $d_{\mathcal{M}_{n}}(\mu,\nu)=2$.
However,
both the strings $1011000000^{n-9}$ and $1101000000^{n-9}$ are not vertices of $\mathcal{M}_{n}$.
So $d_{\mathcal{M}_{n}}(\mu,\nu)>2$,
a contradiction.
Thus,
$\mathcal{M}_{n}$ is not a partial cube for $n\geq9$.
For $n\leq8$,
it is clear that there are at most three bits of $1$s in a vertex of $\mathcal{M}_{n}$.
Suppose $\alpha=a_{1}a_{2}\ldots a_{n},
\beta=b_{1}b_{2}\ldots b_{n}\in V(\mathcal{M}_{n})$,
$H(\alpha,\beta)=t$ and $a_{i_{j}}\neq b_{i_{j}}$ for $j\in1:t$.
We claim that $d_{\mathcal{M}_{n}}(\alpha,\beta)=H(\alpha,\beta)$.
Obviously,
the assertion holds for $H(\alpha,\beta)=1$.
Now we suppose that it holds for vertices with Hamming distance $t-1$ ($t\geq2$).
Without loss of generality,
assume that $a_{i_{1}}=1$ and $b_{i_{1}}=0$.
There are two cases:
at least one of $a_{i_{1}-1}$ and $a_{i_{1}+1}$ is $0$,
both $a_{i_{1}-1}$ and $a_{i_{1}+1}$ are $1$s.
Note that here if $i_{1}=1$,
then $i_{1}-1=n$
and if $i_{1}=n$,
then $i_{1}+1=1$.
For the first case,
let $\alpha'$ be the string obtained from $\alpha$ by changing $a_{i_{1}}$ from $1$ to $0$.
Then $\alpha\in V(\mathcal{M}_{n})$ and $H(\alpha',\beta)=t-1$.
For the second case,
let $\alpha'$ be the string obtained from $\alpha$ by changing $a_{i_{1}+1}$ from $1$ to $0$.
Since $a_{i_{1}+2}$ must be $0$,
$\alpha\in V(\mathcal{M}_{n})$ and $H(\alpha',\beta)=t-1$.
By hypothesis,
$d_{\mathcal{M}_{n}}(\alpha,\beta)
=1+d_{\mathcal{M}_{n}}(\alpha',\beta)
=1+t-1
=t$.
This completes the proof.
$\Box$

It is shown (cf. the proof of Proposition 3.7 of in \cite{HIK}) that the median of the triple  in $Q_{n}$ can be obtained by the \emph{majority rule}:
the $i$-th coordinate of the median is equal to the element that appears at least twice among the $i$-th coordinate of this triple.

\trou \noi {\bf Theorem 5.2.}
\emph{Graph $\mathcal{M}_{n}$ is a median graph if and only if $n\leq7$.}

\trou \noi {\bf Proof.}
For $n\geq7$,
let
$\alpha=11100000^{n-7},
\beta=10000000^{n-7}$
and $\gamma=00100000^{n-7}$ be vertices of $\mathcal{M}_{n}$ embedded into hypercube $Q_{n}$.
Then by majority rule,
the unique candidate for a median is $\eta=10100000^{n-7}$.
Obviously,
$\eta$ is not a vertex of $\mathcal{M}_{n}$,
hence $\mathcal{M}_{n}$ is not a median graph for $n\geq7$.
As $\mathcal{M}_{n}$ is a star for $n\in1:4$,
$\mathcal{M}_{n}$ is median graph obviously.
Now suppose that $n=5$ or $6$.
Then there are at most two bits of $1$ in a vertex of $\mathcal{M}_{n}$.
Let $\alpha=a_{1}a_{2}\ldots a_{n},$
$\beta=b_{1}b_{2}\ldots b_{n},$
$\gamma=c_{1}c_{2}\ldots c_{n}$ be arbitrary vertices of $\mathcal{M}_{n}$ embedded into $Q_{n}$,
and $\omega$ be the string got by majority rule.
Since the $n$ is at most $6$,
the weight of $\omega$ is not more than 2.
If at most one of $i\in1:n$,
the majorities of $a_{i}, b_{i}, c_{i}$ is 1,
then $\omega\in V(\mathcal{M}_{n})$.
Assume that for $i,j\in1:n$,
the majorities of $a_{i}, b_{i}, c_{i}$ and $a_{j}, b_{j}, c_{j}$ are both equal to 1.
If $j=i+1$,
then $\omega\in V(\mathcal{M}_{n})$.
If $j>i+1$,
then the $i$-th and the $j$-th of at least one of $\alpha,\beta,\gamma$ are both $1$.
This means that $j-i>2$ and so $\omega\in V(\mathcal{M}_{n})$.
This completes the proof.
$\Box$

\section{Further directions}

In this section,
some questions and ideas on Fibonacci and Lucas $p$-cubes are proposed for further investigation.

It is shown that every Fibonacci-run graph $\mathcal{R}_{n}$ has a Hamiltonian path,
and conjecture that $\mathcal{R}_{n}$ is Hamiltonian if and only if $n\equiv1$ (mod 3) \cite{EI1}.
On graph $\mathcal{R}_{n}$,
the following question arises.
From the decomposition of graph $\mathcal{M}_{n}$ as shown in Fig. 3,
it can be find that there are many vertices such that with degree 2.
There do not  seem much likelihood of the graph $\mathcal{M}_{n}$ being Hamiltonian.

\trou \noi {\bf Question 6.1.}
\emph{Which Associated Mersenne graphs $\mathcal{M}_{n}$ are Hamiltonian or have a Hamiltonian path $?$}

The nature of the degree sequences of Fibonacci-run graphs $\mathcal{R}_{n}$ was considered \cite{EI2},
and a refinement of the generating function of the degree sequences $\mathcal{R}_{n}$ was obtained.
In relation to this,
a natural question is raised as the following.

\trou \noi {\bf Question 6.2.}
\emph{For given $n$ and $k$,
how many vertices of Associated Mersenne graph $\mathcal{M}_{n}$ have degree $k?$}

For a graph $G$,
let $n\geq0$ and $c_{n}(G)$ be the number of induced subgraphs of $G$ isomorphic to $Q_{n}$.
The cube polynomial $C(G, x)$ of $G$,
is the corresponding counting polynomial,
that is,
the generating function $C(G, x)=\mathop{\sum}\limits_{n\geq0}c_{n}(G)x^{n}$.
It was introduced in \cite{BKS}.
The cube polynomial of Fibonacci and Lucas cubes has been determined \cite{KM2}.
The cube polynomial of Fibonacci and Lucas $p$-cubes was considered in \cite{WeiY}.
The question corresponding to Associated Mersenne graph is as following.

\trou \noi {\bf Question 6.3.}
\emph{What is cube polynomial for Associated Mersenne graph $\mathcal{M}_{n}?$}

\section*{Appendix}

Figures of associated Mersenne graphs $\mathcal{M}_{n}$ for $n\in7:10$ are given in this appendix,
see Figs. 4-7.

\begin{center}
\includegraphics[scale=0.60]{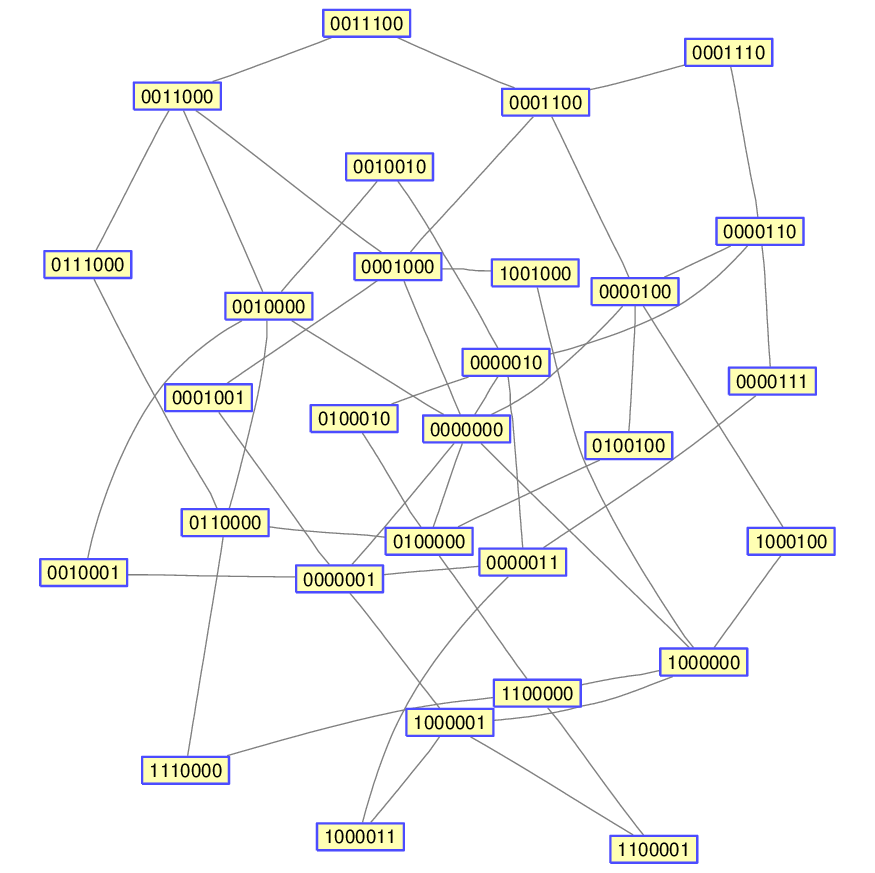}\\
{\footnotesize Fig. 4. The graph $\mathcal{M}_{7}$.}
\end{center}

\begin{center}
\includegraphics[scale=0.55]{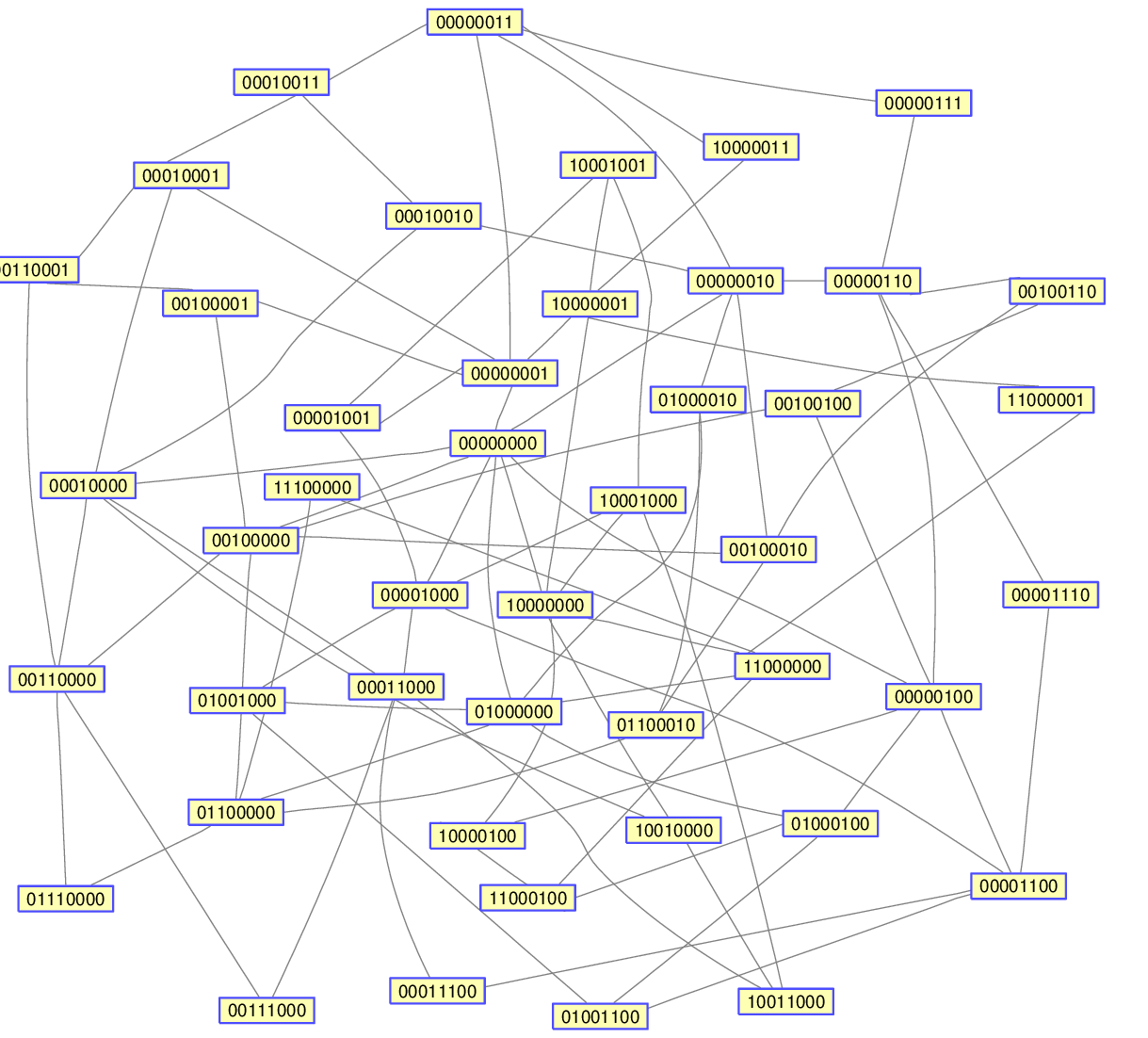}\\
{\footnotesize Fig. 5. The graph $\mathcal{M}_{8}$.}
\end{center}

\begin{center}
\includegraphics[scale=0.55]{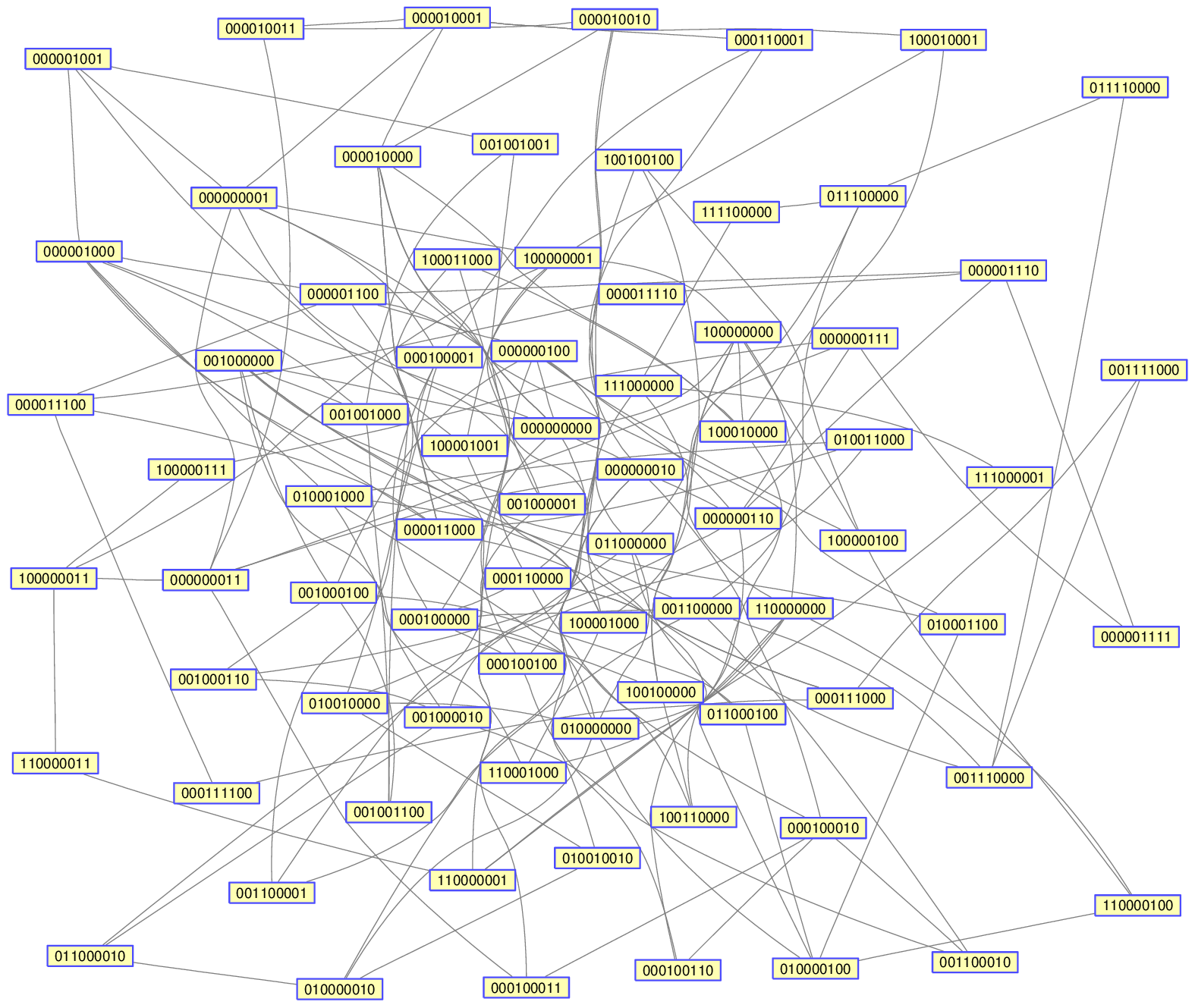}\\
{\footnotesize Fig. 6. The graph $\mathcal{M}_{9}$.}
\end{center}

\begin{center}
\includegraphics[scale=0.50]{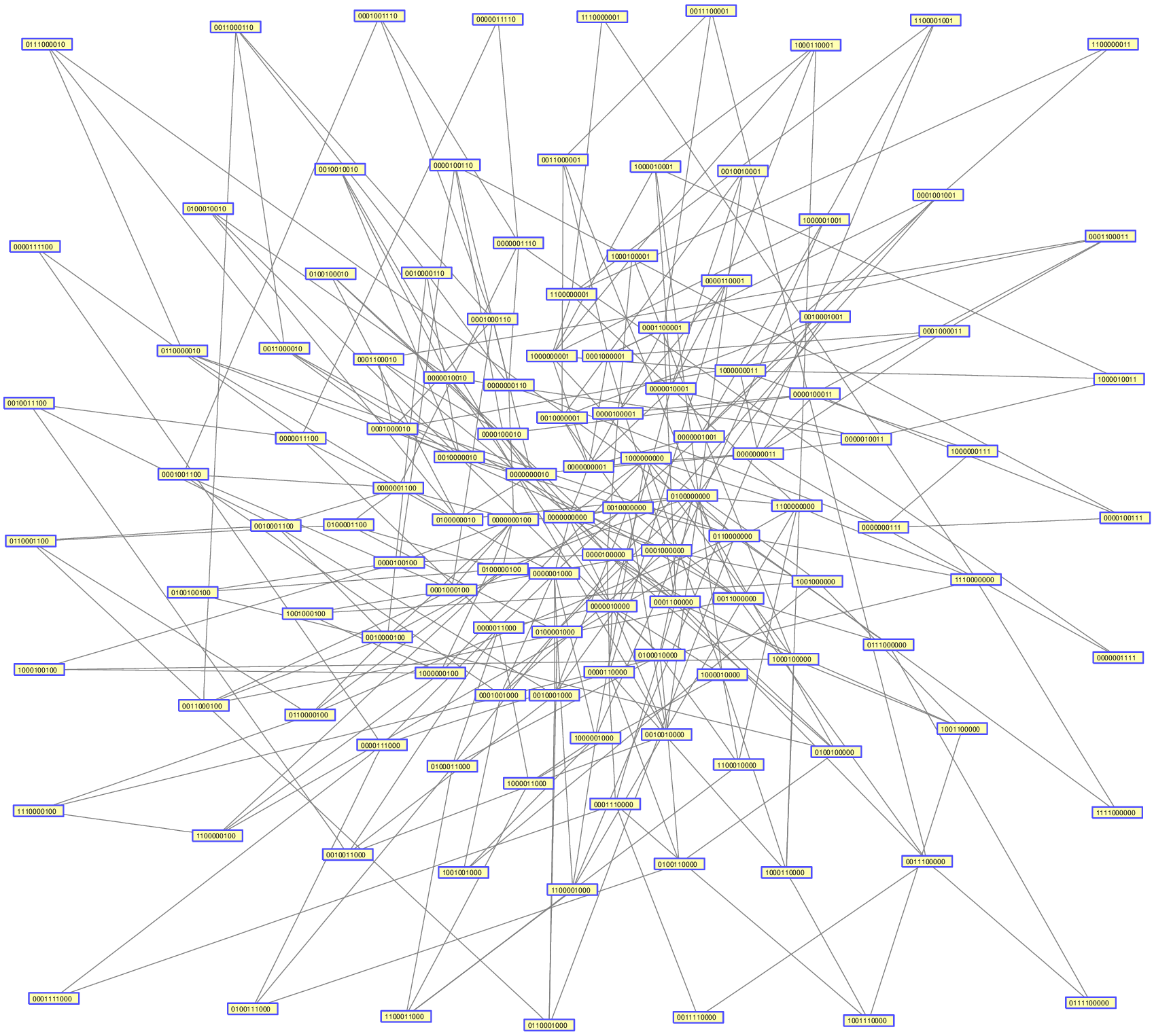}\\
{\footnotesize Fig. 7. The graph $\mathcal{M}_{10}$.}
\end{center}


\begin{thebibliography}{99}

\bibitem{BM} J. A. Bondy, U. S. R. Murty,
Graph Theory with Applications, American Elsevier, New York, 1976.

\bibitem{BKS} B. Bre\v{s}ar, S. Klav\v{z}ar, R. \v{S}krekovski,
The cube polynomial and its derivatives: the case of median graphs,
Electron. J. Comb. 10 (3) (2003), 11 pp. (electronic).

\bibitem{CM}
A. Castro, M. Mollard,
The eccentricity sequences of Fibonacci and Lucas cubes,
Discrete Math. 312 (2012) 1025--1037.

\bibitem{DTZ}
E. Ded\'{o}, D. Torri, N. Zagaglia Salvi,
The observability of the Fibonacci and the Lucas cubes,
Discrete Math. 255 (2002) 55--63.

\bibitem{EI1} {\"O}. E\v{g}ecio\v{g}lu, V. Ir\v{s}i\v{c},
Fibonacci-run graphs I: Basic properties,
Discrete Appl. Math. 295 (2021) 70--84.

\bibitem{EI2} {\"O}. E\v{g}ecio\v{g}lu, V. Ir\v{s}i\v{c},
Fibonacci-run graphs II: Degree sequences,
Discrete Appl. Math. 300 (2021) 56--71.

\bibitem{ESS2} {\"O}. E\v{g}ecio\v{g}lu, E. Sayg{\i}, Z. Sayg{\i},
$k$-Fibonacci cubes: A family of subgraphs of Fibonacci cubes,
Internat. J. Found. Comput. Sci. 31 (5) (2020) 639-661.

\bibitem{ESS1} {\"O}. E\v{g}ecio\v{g}lu, E. Sayg{\i}, Z. Sayg{\i},
Alternate Lucas Cubes,
Internat. J. Found. Comput. Sci. 32 (7) (2021) 871--899.

\bibitem{ESS3} {\"O}. E\v{g}ecio\v{g}lu, E. Sayg{\i}, Z. Sayg{\i},
The structure of $k$-Lucas cubes,
Hacet. J. Math. Stat.
50 (3) (2021) 754--769.

\bibitem{HIK} R. Hammack, W. Imrich, S. Klav\v{z}ar,
Hand book of Product Graphs, second ed., CRC Press, Boca Raton, FL, 2011.

\bibitem{Hsu} W.-J. Hsu,
Fibonacci cubes $-$ a new interconnection topology,
IEEE Trans. Parallel Distrib. Syst. 4 (1993) 3--12.

\bibitem{IKR1} A. Ili\'{c}, S. Klav\v{z}ar, Y. Rho,
Generalized Fibonacci cubes,
Discrete Math. 312 (2012) 2--11.

\bibitem{IKR2} A. Ili\'{c}, S. Klav\v{z}ar, Y. Rho,
Generalized Lucas cubes,
Appl. Anal. Discrete Math. 6 (2012) 82--94.

\bibitem{IM}
A. Ili\'{c}, M. Milo\v{s}evi\'{c},
The parameters of Fibonacci and Lucas cubes,
Ars Math. Contemp. 12 (2017) 25--29.

\bibitem{K1} S. Klav\v{z}ar,
Structure of Fibonacci cubes: a survey,
J. Comb. Optim. 25 (2013) 505--522.

\bibitem{KM} S. Klav\v{z}ar, M. Mollard,
Daisy cubes and distance cube polynomial,
Europ. J. Combin. 80 (2019) 214--223.

\bibitem{KM2}
S. Klav\v{z}ar, M. Mollard,
Cube polynomial of Fibonacci and Lucas cubes,
Acta Appl Math (2012) 117:93--105.


\bibitem{MPS}
E. Munarini, C. Perelli Cippo, N. Zagaglia Salvi,
On the lucas cubes, Fibonacci Quart. 39 (2001) 12--21.

\bibitem{Mu}
E. Munarini,
Pell graphs,
Discrete Math. 342 (2019) 2415--2428.

\bibitem{Sl}
N.J.A. Sloane,
On-Line Encyclopedia of Integer Sequences, http://oeis.org/.

\bibitem{Wei} Jianxin Wei,
The radius and center of Fibonacci-run graphs,
Discrete Applied Mathematics,
325 (2023) 93-96.

\bibitem{WeiY} Jianxin Wei, Yujun Yang,
Fibonacci and Lucas p-cubes,
Discrete Applied Mathematics,
322 (2022) 365-383.

\end{thebibliography}
\end{document}